\documentclass{article}
\usepackage{graphicx} 
\usepackage[margin=1in]{geometry}
\usepackage{amsmath, amsthm, mathtools, amssymb, amsfonts}
\usepackage{tikz}
\usepackage{tikz-cd}
\usepackage{stmaryrd}
\usepackage{bold-extra}
\SetSymbolFont{stmry}{bold}{U}{stmry}{m}{n}

\newcommand\Tstrut{\rule{0pt}{2.6ex}}         
\newcommand\Bstrut{\rule[-0.9ex]{0pt}{0pt}}   
\usepackage{array}
\usetikzlibrary{matrix}
\usepackage{accents}
\usepackage[ordering=Kac,
indefinite-edge={draw=black,fill=white, dotted},
indefinite-edge-ratio=1.5,
mark=o,
root-radius=.1cm]{dynkin-diagrams}
\pgfkeys{/Dynkin diagram,
edge length=.9cm,
fold radius=.6cm,
indefinite edge/.style={
draw=black,
fill=white,
thick, dotted}}
\usepackage[colorinlistoftodos,textsize=small]{todonotes}
\newcommand{\spc}{\hspace{0.1em}}
\newcommand{\bigplus}{\raisebox{0.7ex}{\scalebox{1.1}{$\scriptstyle\boldsymbol{+}$}}}
\newcommand{\bigminus}{\raisebox{0.7ex}{\scalebox{1.1}{$\scriptstyle\boldsymbol{-}$}}}

\makeatletter
\newcommand*\bigcdot{\mathpalette\bigcdot@{.5}}
\newcommand*\bigcdot@[2]{\mathbin{\vcenter{\hbox{\scalebox{#2}{$\m@th#1\bullet$}}}}}
\makeatother

\DeclareUnicodeCharacter{2212}{-}

\theoremstyle{plain}
\newtheorem{theorem}{Theorem}[section]

\newtheorem{lemma}[theorem]{Lemma}
\newtheorem{proposition}[theorem]{Proposition}

\theoremstyle{definition}
\newtheorem{definition}[theorem]{Definition}
\newtheorem{remark}[theorem]{Remark}
\newtheorem{example}[theorem]{Example}

\usepackage{pgfgantt}
\usepackage[
backend=biber,
style=alphabetic,
sorting=nyt
]{biblatex}
\usepackage{lipsum}
\newcommand\blfootnote[1]{%
  \begingroup
  \renewcommand\thefootnote{}\footnote{#1}%
  \addtocounter{footnote}{-1}%
  \endgroup
}

\begin{document}
\title{Odd magical triples and maximal Higgs bundles}
\author{Enya Hsiao}
\date{}

\maketitle

\begin{abstract}
    We introduce the notion of extended magical $\mathfrak{sl}_2$-triples, which is a generalization of magical $\mathfrak{sl}_2$-triples in \cite{Bradlow_Collier_García-Prada_Gothen_Oliveira_2024} and show that, apart from the magical triples in \cite{Bradlow_Collier_García-Prada_Gothen_Oliveira_2024}, which are all even, there are precisely three odd triples of nontube type Hermitian Lie algebras that are extended magical.
    We then show that the intersection of the Slodowy slice of an odd extended magical triple of $G^\textbf{R}$ with the $G^\textbf{R}$-Higgs bundle moduli space yields precisely the maximal components. Finally, assuming that the underlying curve has sufficiently large genus, we give a geometric characterization of extended magical triples and prove a Cayley correspondence for the maximal components of $G^\textbf{R}$-Higgs bundles for nontube type Hermitian Lie groups. \blfootnote{\textit{Date:} January 24, 2026. \\
    \indent \,\,\, \textit{Affliation:} \textsc{Max-Planck-Institute Mathematics in the Sciences, 04103 Leipzig, Germany} \\
    \indent \,\,\, \textit{Email address:} \texttt{enya.hsiao@mis.mpg.de}}
\end{abstract} 
   
\setcounter{tocdepth}{1}
\tableofcontents

\makeatletter
\providecommand\@dotsep{5}
\def\listtodoname{TO-Dos}
\def\listoftodos{\@starttoc{tdo}\listtodoname}
\makeatother

\section{Introduction}
Higher Teichm\"uller theory is the study of connected components of the character variety $$\chi(G^\textbf{R}):=\mathrm{Hom}(\pi_1(X),G^\textbf{R})\sslash G^\textbf{R}$$ of surface group representations into semisimple real Lie groups of \textit{higher} rank $
>1$ that consists entirely of discrete and faithful representations. It emerged in the last decade as a new field in mathematics and is quickly gaining traction and relevance due to an abundance of recent fruitful results that interconnect a significant number of different mathematical fields. 

The first two encounters of higher Teichm\"uller components, and indeed the cornerstone of many examples that followed, were Hitchin's discovery of the Hitchin components for split real groups \cite{Hitchin_1992}, and Burger-Iozzi-Wienhard's seminal work on the components of maximal representations for Lie groups of Hermitian type \cite{Burger_Iozzi_Wienhard_2010}. They were shown to consist entirely of discrete and faithful representations using the notion of Anosov representations introduced by Labourie \cite{Labourie_2006}. Subsequent work by Labourie, Guichard and Wienhard \cite{guichard2025positivityrepresentationssurfacegroups}\cite{Guichard_Wienhard_2025} lead to the introduction of $\Theta$-positive structures on semisimple real Lie groups and that of $\Theta$-positive representations, giving rise to higher Teichm\"uller components \cite{beyrer2024positivitycrossratioscollarlemma}. 

Higgs bundles enter the picture as an important tool in higher Teichm\"uller theory via the nonabelian Hodge correspondence, a homeomorphism between the character variety and the moduli space of Higgs bundles. Higgs bundle methods were successfully used to count and identify higher Teichm\"uller components \cite{Aparicio-Arroyo_Bradlow_Collier_García-Prada_Gothen_Oliveira_2019},\cite{Bradlow_García-Prada_Gothen_2007},\cite{Collier_2017},\cite{Gothen_1995}. The culmination of these efforts lead to the introduction and classification of magical $\mathfrak{sl}_2$-triples, where an $\mathfrak{sl}_2$-triple $\rho$ of a complex simple Lie algebra $\mathfrak{g}$ is called \textit{magical} if its magical involution $\sigma_\rho:\mathfrak{g}\to \mathfrak{g}$ given by \eqref{eqn:magical_involution} is a Lie algebra involution. An $\mathfrak{sl}_2$-triple $\hat{\rho}$ of a real form $\mathfrak{g}^\textbf{R}\subset \mathfrak{g}$ is called magical if its Cayley transform $\rho:\mathfrak{sl}_2\textbf{C} \to \mathfrak{g}$ is magical and $\mathfrak{g}^\textbf{R}$ is the canonical real form fixed by the Lie algebra involution $\sigma_\rho$. 

Remarkably, given a magical triple $\hat{\rho}$ of $G^\textbf{R} \subset G$, the real form of a complex simple Lie group for which the magical involution integrates to an involution on the group level, there is an injective, open and closed map on the level of $K^p$-twisted moduli spaces \cite[Theorem 7.1]{Bradlow_Collier_García-Prada_Gothen_Oliveira_2024} whose image is a union of higher Teichm\"uller components:
\begin{align}
\label{eqn:Cayley_correspondence}   \Psi_\rho:\mathcal{M}_{K^{m_c+1}}(\tilde{G}^\textbf{R})\times \prod_{j=1}^{r(\rho)} \mathcal{M}_{K^{l_j+1}}(\textbf{R}^+) \to \mathcal{M}(G^\textbf{R})
\end{align}
where $\tilde{G}^\textbf{R}$ is the semisimple part of a certain subgroup of $G^\textbf{R}$, $r(\rho)=\mathrm{rank}(\mathfrak{g}(\rho))$ is the rank of a certain subalgebra of $\mathfrak{g}$ and the integers $m_c$ and $l_j$ are weights associated to the $\mathfrak{sl}_2$-module decomposition $\mathfrak{g}=W_{m_c}\bigoplus_{l_j} W_{l_j}$. This phenomenon is called the \textit{generalized Cayley correspondence} and the connected components in the image of $\Psi_\rho$ are called a Cayley components. 

In \cite{Bradlow_Collier_García-Prada_Gothen_Oliveira_2024} it is conjectured that the Cayley components are, apart from the maximal components for nontube type Hermitian Lie groups, precisely all of the higher Teichm\"uller components. In this article, we prove that the maximal components for nontube type Hermitian Lie groups are also described by a generalized Cayley correspondence by extending the notion of magical $\mathfrak{sl}_2$-triples.

Very briefly, we say that an $\mathfrak{sl}_2$-triple is an \textit{extended magical triple} if there exists a Lie algebra involution $\sigma$ that agrees with the magical involution on the even weight spaces of the $\mathfrak{sl}_2$-module decomposition $\mathfrak{g}=\bigoplus_{j} W_j$. An $\mathfrak{sl}_2$-triple $\hat{\rho}$ of a real form $\mathfrak{g}^\textbf{R} \subset \mathfrak{g}$ will be called \textit{extended magical} if its Cayley transform $\rho:\mathfrak{sl}_2\textbf{C}\to \mathfrak{g}$ is extended magical and $\mathfrak{g}^\textbf{R}$ is canonical real form fixed by $\sigma$.

With this definition, the magical triples in the sense of \cite{Bradlow_Collier_García-Prada_Gothen_Oliveira_2024} are naturally extended magical and will be called \textit{even magical triples}, since they all turn out to be even $\mathfrak{sl}_2$-triples. Those that are not captured by the original definitions are automatically odd and are called \textit{odd magical triples}. Our first main result is the classification of real forms $\mathfrak{g}^\textbf{R}$ admitting odd magical triples:

\begin{theorem}[Theorem \ref{thm:classification_theorem}]
An $\mathfrak{sl}_2$-triple $\hat{\rho}:\mathfrak{sl}_2\to \mathfrak{g}^\textbf{R}$ is an odd magical triple if and only if $\mathfrak{g}^\textbf{R}$ is one of the following Hermitian Lie algebra of nontube type:
\begin{center}
        \begin{tabular}{|c|c|c|c|}
        \hline $\mathfrak{g}$ & $\mathfrak{sl}_{p+q}\mathrm{\textbf{C}}$ & $\mathfrak{so}_{4m+2}\mathrm{\textbf{C}}$ & $E_6$ \Tstrut  \\ [0.2ex]
                 \hline 
$\mathfrak{g}^\textbf{R}$ & $\mathfrak{su}(p,q)$ & $\mathfrak{so}^*_{4m+2}$ & $E_6^{-14}$ \Tstrut \\ [0.2ex] \hline
    \end{tabular}
\end{center}
\end{theorem}
\medskip
 
Let $\mathcal{E}_{G^\textbf{R}}=(E_H,\varphi)$ be a $G^\textbf{R}$-Higgs bundle of a Hermitian Lie group with maximal compact subgroup $H^\textbf{R}$. The \textit{Toledo invariant} of $\mathcal{E}_{G^\textbf{R}}$ is defined to be the degree 
$\tau=\deg(E_H(\tilde{\chi}))$ of the line bundle associated to the \textit{lifted Toledo character} $\tilde{\chi}:H\to \textbf{C}$ (see Definition \ref{def:toledo_character}), and is an integer satisfying the Milnor--Wood inequality: $|\tau|\le \mathrm{rk}(G^\textbf{R}/H^\textbf{R})(2g-2)$. A $G^\textbf{R}$-Higgs bundle with maximal Toledo number is called \textit{maximal} and a connected component consisting of maximal $G^\textbf{R}$-Higgs bundles is called a \textit{maximal component}. The maximal components $\mathcal{M}_{\text{max}}(G^\textbf{R})\subset \mathcal{M}(G^\textbf{R})$ are known to be higher Teichm\"uller components in the $G^\textbf{R}$-character variety \cite{Burger_Iozzi_Wienhard_2010}.

On the other hand, for any $\mathfrak{sl}_2$-triple $\rho$ of $G$, let $C$ be the centralizer of the image $\rho(\mathrm{SL}_2)\subset G$ and $V_\rho^{>0}:=\bigoplus_{j>0} V_j$ be the sum of the nontrivial highest weight spaces $V_j\subset W_j$. The \textit{$\rho$-Slodowy category} $\mathcal{B}_\rho(G)$ introduced by Collier and Sanders \cite{Collier_Sanders_2021} consists of objects which are tuples $((\mathcal{E}_C,\varphi_C),\{\varphi_{n_j}\})$, where $(\mathcal{E}_C,\varphi_C)$ is a $C$-Higgs bundle and $\varphi_{j}\in H^0(X, \mathcal{E}_C[V_j]\otimes K)$ are holomorphic sections taking values in $V_j$. The Slodowy map $\hat{\Psi}_\rho$ is an assignment that sends each tuple in the Slodowy category to a $G$-Higgs bundle 
\begin{align*}
\hat{\Psi}_\rho:\mathrm{Slo}_\rho \to \mathcal{H}(G), \quad 
    ((E_C,\varphi_C),\{\varphi_j\}) \mapsto ((E_C\times E_T)[G], \underline{f}+\varphi_C+\sum_{j=1}^N \varphi_{n_j}),
\end{align*}
where $(E_T,\underline{f})$ is the uniformizing bundle over $X$ associated to $\rho$.
For a real form $G^\textbf{R}$, we will call the image of the Slodowy map restricted to the preimage of polystable $G^\textbf{R}$-Higgs bundles in $\mathcal{M}(G^\textbf{R})$ the \textit{Slodowy slice} $\mathrm{Slo}_\rho$ associated to $\rho$. Our next result shows that the Slodowy slice of an odd magical triple correspond to the maximal components:

\begin{theorem}[Theorem \ref{thm:maximal_components}]
    Let $G$ be a simple complex Lie group and let $\hat{\rho}$ be an odd magical triple of the real form $G^\textbf{R}\subset G$ with $\rho$ its Cayley transform. The Slodowy slice $\mathrm{Slo}_\rho \subset \mathcal{M}(G^\textbf{R})$ is precisely the union of maximal $G^\textbf{R}$-Higgs bundles $\mathcal{M}_{\mathrm{max}}(G^\textbf{R})$.
\end{theorem}

Combining the above with a result in \cite{Kim_Pansu_2014} stating that given sufficiently large genus of the underlying curve $X$, all of the connected components of the $G^\textbf{R}$-character variety achieves expected dimension except for the maximal components when $G^\textbf{R}$ is Hermitian of nontube type, we obtain a geometric characterization of extended magical triples:

\begin{theorem}[Theorem \ref{thm:characterization}]
Let $G^\textbf{R} \subset G$ be a real form of a complex simple Lie group with Lie algebra $\mathfrak{g}^\textbf{R}$. If an $\mathfrak{sl}_2$-triple $\hat{\rho}$ of $\mathfrak{g}^\textbf{R}$ is extended magical, then its Slodowy slice $\mathrm{Slo}_\rho \subset \mathcal{M}(G^\textbf{R})$ is a union of connected components. Conversely, assuming that $X$ has genus $g(X)\ge 2\dim_\textbf{R}(G^\textbf{R})^2$, then $\mathrm{Slo}_\rho \subset \mathcal{M}(G^\textbf{R})$ is a union of connected components only if the $\mathfrak{sl}_2$-triple $\hat{\rho}$ of $\mathfrak{g}^\textbf{R}$ is extended magical.
\end{theorem}

Finally, given an odd magical triple $\hat{\rho}$ of $G^\textbf{R}$, we prove a Cayley correspondence for the maximal components of $G^\textbf{R}$-Higgs bundles when $G^\textbf{R}$ is Hermitian of nontube type:
\begin{theorem}[Theorem \ref{thm:cayley_correspondence_odd}]
    Let $G$ a complex simple Lie group and let $\hat{\rho}$ be an odd magical triple of the real form $G^\textbf{R}\subset G$. Let $\tilde{G}^\textbf{R}$ be the subgroup of $G^\textbf{R}$ with Lie algebra $\tilde{\mathfrak{g}}^\textbf{R} \subset \mathfrak{g}^\textbf{R}$ and maximal compact subgroup $\mathcal{C}:=C\cap G^\textbf{R}$. Let $\hat{\rho}$ be an odd magical triple of $G^\textbf{R}$. Consider $\tilde{\mathfrak{g}}$, the semisimple part of $\mathfrak{g}_0$ and let $\tilde{\mathfrak{g}}^\textbf{R}\subset \tilde{\mathfrak{g}}$ be the real form defined with respect to the involution \eqref{eq:theta_involution} restricted to $\tilde{\mathfrak{g}}$. Let $\tilde{G}^\textbf{R}$ be the subgroup of $G^\textbf{R}$ with Lie algebra $\tilde{\mathfrak{g}}^\textbf{R}$ and maximal compact subgroup $C\cap G^\textbf{R}$. There is an injective, open and closed map on the level of $K^p$-twisted Higgs bundle moduli spaces:
    \begin{align*}
\Psi_\rho:\mathcal{M}_{K^2}(\tilde{G}^\textbf{R})\times  \mathcal{M}(\textbf{R}^+) \to \mathcal{M}_{\mathrm{max}}(G^\textbf{R})\subset \mathcal{M}(G^\textbf{R})
    \end{align*}
    whose image is the union of maximal components.
\end{theorem}

The contents of this paper is organized as follows: In Section 2 we set up some notation and briefly review the necessary background for the subsequent sections. In Section 3 we prove the classification of extended magical triples and in Section 4 we exhibit explicitly the $\mathfrak{sl}_2$-data associated to the odd magical triples. In Sections 5 and 6, we prove a geometric characterization of extended magical triples and a Cayley correspondence for the maximal components of $G^\textbf{R}$-Higgs bundles for nontube type Hermitian Lie groups.
\subsection*{Acknowledgement}
This work was inspired by questions posed by Colin Davalo and Mengxue Yang during a research visit at the Kavli institute in February 2025 and the tri-semester at the Institut Henri Poincar\'e in April 2025. The author is grateful for their ideas and stimulating discussion and gives many thanks to Brian Collier and Junming Zhang for useful remarks. The author is also tremendously indebted to her PhD advisor Anna Wienhard for helpful insights and encouragement, and to her mentor Georgios Kydonakis for carefully reading through the manuscript, providing corrections and feedback. The author is funded by the Deutsche 
Forschungsgemeinschaft (DFG, German Research
Foundation) – 541679129. This project has received funding from the European Research Council (ERC) under the European Union’s Horizon 2020 research and innovation programme (grant agreement No 101018839).

\section{Preliminaries}
\subsection{Notations and conventions}
\label{sec:notations}
As a conventional rule, we will denote complex Lie algebras by unadorned fraktur letters such as $\mathfrak{g}$, $\mathfrak{h}$ and $\mathfrak{m}$ and real Lie algebras with boldface superscript such as $\mathfrak{g}^\textbf{R}$, $\mathfrak{h}^\textbf{R}$ and $\mathfrak{m}^\textbf{R}$. Similarly, real Lie groups will be indicated by a boldfaced superscript (e.g. $G^\textbf{R}$), while complex Lie groups will not. As a rule of thumb, the unadorned complex object, e.g. $\mathfrak{g}$ or $G$, will be the complexification of real object of the same letter, e.g. $\mathfrak{g}^\textbf{R}$ or $G^\textbf{R}$. Complex $\mathfrak{sl}_2$-triples will be denoted by $\rho:\mathfrak{sl}_2\mathbf{C}\to \langle e,h,f\rangle\subset \mathfrak{g}$, while real $\mathfrak{sl}_2$-triples will be decorated with a hat $\hat{\rho}:\mathfrak{sl}_2\mathrm{\textbf{R}}\to \langle \hat{e},\hat{h},\hat{f} \rangle \subset\mathfrak{g}^\textbf{R}$. Whenever we are dealing with a Hermitian Lie group or algebra of nontube type, we will fix a maximal subtube and denote it and its subgroups or subalgebra with the subscript $t$, e.g. $G_t^\textbf{R}$, $H^\textbf{R}_t$ or $\mathfrak{g}_t^\textbf{R}$, $\mathfrak{h}^\textbf{R}_t$. For statements concerning extended magical triples and the restricted Slodowy map, whenever we choose an ambient Lie group $G$ whose Lie algebra is $\mathfrak{g}$, we ask that $G$ be connected, of finite center and such that the magical involution integrates to a Lie group involution. 

Throughout this article, we fix a connected compact Riemann surface $X$ with genus $g\ge 2$ as our base space, and all geometric objects we consider will be over $X$. Let $K:=T^*X$ denote the canonical line bundle, and fix a choice of square root $K^{1/2}$. For an arbitrary group $G$, real or complex, we will denote by $E_G$ a holomorphic principal $G$-bundle and by $\mathcal{E}_G$ a $G$-Higgs bundle over $X$. We will write $E_G(\chi)$ for the line bundle associated to the character $\chi:G\to \textbf{C}$. Given a group homomorphism $G\to G'$ or representation $G\to \mathrm{GL}(V)$, we will write the associated $G'$-bundle or vector bundle as $E_G[G']$ or $E_G[V]$ and omit the homomorphism or representation when the context is clear. Given a line bundle $L$, we write $\mathcal{H}_L(G)$ for the configuration space of $L$-twisted $G$-Higgs pairs, $\mathcal{M}_L(G)$ for the moduli space of $L$-twisted polystable $G$-Higgs bundles. When $L=K$, we omit the subscript.

\subsection{$\mathfrak{sl}_2$-triples}
We briefly review here the necessary background pertaining to $\mathfrak{sl}_2$-triples that we require later on, and refer the reader to \cite[Section 2]{Bradlow_Collier_García-Prada_Gothen_Oliveira_2024} and \cite{Collingwood_McGovern_1993} for more detailed accounts. 

\subsubsection{Complex $\mathfrak{sl}_2$-triples}
\label{sec:sl2}
Let $\mathfrak{g}$ be a simple complex Lie algebra. An \textit{$\mathfrak{sl}_2$-triple} of $\mathfrak{g}$ is an embedding $\rho:\mathfrak{sl}_2\textbf{C} \to \langle e,h,f\rangle \subset \mathfrak{g}$, where the generators satisfy the relations
\begin{align*}
[h,e]=2e, \,\, [h,f]=-2f, \,\, [e,f]=h. \end{align*}
Unlike the set up in \cite{Bradlow_Collier_García-Prada_Gothen_Oliveira_2024}, we do not exclude the trivial $\mathfrak{sl}_2$-triple from our considerations.

Let $G$ be a connected simple Lie group with Lie algebra $\mathfrak{g}$. The Jacobson--Morozov Theorem establishes a one-to-one correspondence between $G$-conjugacy classes of $\mathfrak{sl}_2$-triples and nilpotent orbits
\begin{align*}\{\rho:\mathfrak{sl}_2\textbf{C} \hookrightarrow \mathfrak{g}\}/G
     \,\,\xleftrightarrow{1:1} \,\, \{e \in \mathfrak{g} \,\,\text{nilpotent}\}/G.
\end{align*}
Therefore, up to $G$-conjugacy an $\mathfrak{sl}_2$-triple $\rho$ is determined uniquely by a nilpotent element $e \in \mathfrak{g}$ and vice versa. \\

An $\mathfrak{sl}_2$-triple $\rho:\mathfrak{sl}_2\textbf{C}\to \langle e,h,f\rangle \subset \mathfrak{g}$ induces two compatible decompositions of $\mathfrak{g}$ that we now describe. On the one hand, the semisimple element $h$ gives rise to a $\mathbf{Z}$-grading on $\mathfrak{g}$ 
$$\mathfrak{g}=\bigoplus_{j=-l}^l \mathfrak{g}_j,$$ 
where $\mathfrak{g}_j$ is the $j$-weight space of $\mathrm{ad}_h$ and we have that $[\mathfrak{g}_j,\mathfrak{g}_i]\subset \mathfrak{g}_{i+j}$. On the other hand, the Lie algebra $\mathfrak{g}$ is an $\mathfrak{sl}_2\textbf{C}$-representation $\mathfrak{sl}_2\textbf{C}\to \mathfrak{gl}(\mathfrak{g}), \,\, X\mapsto \mathrm{ad}_{\rho(X)}$ induced by the embedding $\rho$, hence admits an $\mathfrak{sl}_2\textbf{C}$-module decomposition into sums of irreducible $\mathfrak{sl}_2\textbf{C}$-representations. We denote by $w_j$ the unique $(j+1)$-dimensional irreducible $\mathfrak{sl}_2\textbf{C}$-representation and let $W_j$ be $n_j$-copies of $w_j$ appearing as summands in the $\mathfrak{sl}_2\mathbf{C}$-module decomposition of $\mathfrak{g}$,
$$\mathfrak{g}=\bigoplus_{(j,n_j)} W_j=\bigoplus_j w_j^{\oplus n_j}.$$
This decomposition is fully characterized by the number of copies $(n_0,n_1,...,n_m)$ of each $w_j$. Let $v_j$ be the highest weight space of $w_j$ and $V_j={v_j}^{\oplus n_j}$ the sum of the highest weight spaces in $W_j$. We denote by $V_\rho=\bigoplus_{j} V_j$ the sum of the highest weight spaces and $V_\rho^{>0}=\bigoplus_{j>0} V_j$ the sum of the \textit{nontrivial} highest weight spaces. 

Now each $w_j$ has a basis of $(j+1)$-vectors, one in each weight space $\mathfrak{g}_{j-2k}$ for $0\le k\le j$. Up to a scaling factor the basis vectors can be risen or lowered with respect to the $\textbf{Z}$-grading via $\mathrm{ad}_f:\mathfrak{g}_i\to \mathfrak{g}_{i-2}$ and $\mathrm{ad}_e:\mathfrak{g}_i\to \mathfrak{g}_{i+2}$. Therefore, the $\textbf{Z}$-grading and $\mathfrak{sl}_2\textbf{C}$-module decompositions of $\mathfrak{g}$ can be combined in the following way:
    \begin{align}
    \label{eqn:sl2_decomposition}
\mathfrak{g} &=\bigoplus_{(j,n_j)} \bigoplus_{k=0}^j W_j \cap \mathfrak{g}_{j-2k} \\ &=\bigoplus_{(j,n_j)} \bigoplus_{k=0}^j \,(\mathrm{ad}_f)^k (V_j),
    \end{align}
where the summands $W_j \cap \mathfrak{g}_{j-2k}=(\mathrm{ad}_f)^k (V_j)$ are $n_j$-dimensional.

To summarize, given an $\mathfrak{sl}_2$-triple $\rho:\mathfrak{sl}_2\textbf{C}\to \mathfrak{g}$, we obtain a set of highest weight spaces $V_j$ and their dimensions $\dim V_j=n_j$, of which $V_0$ coincides with the centralizer $\mathfrak{c}$, the subalgebra of all elements commuting with $\langle e, h, f\rangle$. Furthermore, associated to the $\mathfrak{sl}_2\textbf{C}$-module decomposition is a parabolic subalgebra $\mathfrak{p}$, its opposite parabolic subalgebra $\mathfrak{p}^{\mathrm{opp}}$, a unipotent subalgebra $\mathfrak{u}$, and a Levi subalgebra $\mathfrak{l}$,
\begin{align*}
    \mathfrak{p}=\bigoplus_{i\ge 0} \mathfrak{g}_i, \quad
    \mathfrak{p}^{\text{opp}}=\bigoplus_{i\le 0} \mathfrak{g}_i, \quad\mathfrak{u}=\bigoplus_{i>0} \mathfrak{g}_i, \quad
    \mathfrak{l}=\mathfrak{g}_0.
\end{align*}
All of this information can be completed encoded in the \textit{weighted Dynkin diagram} of $\rho$, which is the Dynkin diagram of $\mathfrak{g}$ with each simple root $\alpha$ labeled by the numbers $\{0,1,2\}$ corresponding to the value $\mathrm{ad}_h(\alpha)$. The weighted Dynkin diagrams is a complete invariant of nilpotent orbits in $\mathfrak{g}$, and therefore of $\mathfrak{sl}_2$-triples up to $G$-conjugacy by the Jacobson-Morozov Theorem. We will call this the \textit{$\mathfrak{sl}_2$-data associated to $\rho$}.

\subsubsection{Real $\mathfrak{sl}_2$-triples}
Let $G^\textbf{R}\subset G$ be a real form with maximal compact subgroup $H^\textbf{R}$ and Cartan decomposition $\mathfrak{g}^\textbf{R}=\mathfrak{h}^\textbf{R}+\mathfrak{m}^\textbf{R}$. Let $H\subset G$ be complexification of $H^\textbf{R}$ with complexified Cartan decomposition $\mathfrak{g}=\mathfrak{h}+\mathfrak{m}$.  

The objects we are interested in this paper are real $\mathfrak{sl}_2$-triples $\hat{\rho}:\mathfrak{sl}_2\mathrm{\textbf{R}}\to \langle \hat{e},\hat{h},\hat{f} \rangle \subset\mathfrak{g}^\textbf{R}$ up to $G^\textbf{R}$-conjugacy. The Jacobson--Morozov Theorem holds over real Lie groups, hence these are in bijection with real nilpotent orbits in $\mathfrak{g}^\textbf{R}$. Now $G^\textbf{R}$-nilpotent orbits in $\mathfrak{g}^\textbf{R}$ are in one-to-one correspondence with $H$-nilpotent orbits in $\mathfrak{m}$ via the \textit{Sekiguchi correspondence}
\begin{align*} \{\hat{e} \in \mathfrak{g}^\textbf{R} \,\,\text{nilpotent}\}/G^\textbf{R}
\,\,\xleftrightarrow{1:1} \,\,
\{e \in \mathfrak{m} \,\,\text{nilpotent}\}/H,
\end{align*}
from which we can study $G^\textbf{R}$-conjugacy classes of real $\mathfrak{sl}_2$-triples $\hat{\rho}$ equivalently as $H$-conjugacy classes of complex $\mathfrak{sl}_2$-triples $\rho$ with $e\in \mathfrak{m}$. The complex triple $\rho$ under this correspondence is called the \textit{Cayley transform} of $\hat{\rho}$.

To be more precise, let $\sigma:\mathfrak{g}\to \mathfrak{g}$ be the Cartan involution determined by $\mathfrak{g}=\mathfrak{h}+\mathfrak{m}$. A real $\mathfrak{sl}_2$-triple $\hat{\rho}$ of $\mathfrak{g}^\textbf{R}$ is called a \textit{Cayley triple} if it satisfies
\begin{align*}
    \sigma(\hat{f})=-\hat{e}, \, \,\sigma(\hat{h})=-\hat{h},\,\,\sigma(\hat{e})=-\hat{f},
\end{align*}
while a complex $\mathfrak{sl}_2$-triple $\rho$ of $\mathfrak{g}$ is called a \textit{normal triple} if 
\begin{align*}
    \sigma(f)=-f, \, \,\sigma(h)=h,\,\,\sigma(e)=-e.
\end{align*}
Notice that the nilpotent element $e$ of a normal triple belongs to $\mathfrak{m}$, the $-1$-eigenspace of $\sigma$. 
 
Any real triple $\hat{\rho}$ is a Cayley triple up to $G^\textbf{R}$-conjugation \cite[Theorem 9.4.1]{Collingwood_McGovern_1993}. Therefore, without loss of generality, we will assume from now on that any real $\mathfrak{sl}_2$-triple $\hat{\rho}$ is a Cayley triple. The Sekiguchi correspondence is shown by exhibiting an explicit bijection mapping $G^\textbf{R}$-conjugacy classes of Cayley triples $\hat{\rho}$ to $H$-conjugacy classes of normal triples $\rho$. Given a real $\mathfrak{sl}_2$-triple $\hat{\rho}$, we will denote by $\rho$ its Cayley transform, which is always a normal triple under the Sekiguchi correspondence. 

\subsection{Maximal Higgs bundles}
One of our main results concerns maximal $G^\textbf{R}$-Higgs bundles when $G^\textbf{R}$ is a simple Hermitian Lie group of nontube type. We supply here the necessary preliminaries, mostly based on \cite{Biquard_García-Prada_Rubio_2017}.

\subsubsection{Hermitian Lie groups $G^\textbf{R}$ of nontube type}

We study connected non-compact simple real Lie groups $G^\textbf{R}$ with finite center and maximal compact subgroup $H^\textbf{R}$ whose associated symmetric space $G^\textbf{R}/H^\textbf{R}$ is an irreducible Hermitian symmetric space. In this case, the center $\mathfrak{z}^\textbf{R}$ of $\mathfrak{h}^\textbf{R}$ is one dimensional, and a special element $J \in \mathfrak{z}^\textbf{R}$ defines an almost complex structure on $G^\textbf{R}/H^\textbf{R}$ compatible with the Riemannian metric, turning it into a K\"ahler manifold. Such a Lie group is called a \textit{Hermitian Lie group}. Accordingly, a real Lie algebra $\mathfrak{g}^\textbf{R}$ is called a \textit{Hermitian Lie algebra} if it is the Lie algebra of a Hermitian Lie group. 

A Hermitian symmetric space $G^\textbf{R}/H^\textbf{R}$ is said to be of tube type if it is biholomorphic to a tube type domain $V+i\Omega$, where $V$ is a real vector space and $\Omega\subset V$ is an open cone. It is said to be of nontube type otherwise. A Hermitian Lie group $G^\textbf{R}$ is of \textit{tube type} (resp. \textit{nontube type}) if its associated symmetric space is of tube type (resp. nontube type). Every Hermitian $G^\textbf{R}$ of nontube type contains a connected subgroup $G_t^\textbf{R} \subset G^\textbf{R}$ called a \textit{maximal subtube}, whose associated symmetric space $G_t^\textbf{R}/H_t^\textbf{R}$ is irreducible Hermitian of tube type and which is maximal in the sense that it is not contained in any other connected subgroup with this property. Any two maximal subtubes are conjugate to each other. 

For our convenience, we will fix a choice of $G_t^\textbf{R}$ and refer to it as \textit{the} maximal subtube of $G^\textbf{R}$. However, none of the proofs and results that follow depend implicitly or explicitly on this choice. In the table below, we list all the irreducible Hermitian symmetric spaces of nontube type with our choice of maximal subtube $G_t^\textbf{R}$ along with its maximal compact subgroup $H_t^\textbf{R}$. 
\begin{table}[htb]
\centering
\begin{tabular}{ |c|c|c|c| } 
\hline $G^\textbf{R}$ & $H^\textbf{R}$ & $G^\textbf{R}_t$ & $H^\textbf{R}_t$ \Tstrut\Bstrut\\
  \hline \hline $\mathrm{SU}(p,q)$ & $\mathrm{S}(\mathrm{U}(p)\times \mathrm{U}(q))$ & $\mathrm{SU}(p,p)$ & $\mathrm{S}(\mathrm{U}(p)\times \mathrm{U}(p))$  \Tstrut \Bstrut\\
\hline $\mathrm{SO}^*(4m+2)$ & $ \mathrm{U}(2m+1)$ & $\mathrm{SO}^*(4m)$ & $\mathrm{U}(2m)$ \Tstrut \Bstrut \\
\hline 
$\mathrm{E}_6^{-14}$ & $\mathrm{Spin}(10)\times_{\mathrm{\textbf{Z}}_4} U(1)$ & $\mathrm{Spin}_0(2,8)$ & $\mathrm{Spin}(2)\times \mathrm{Spin}(8)/\mathrm{\textbf{Z}}_2$ \Tstrut\Bstrut \\
\hline
\end{tabular}
\caption{Irreducible Hermitian symmetric spaces $G^\textbf{R}/H^\textbf{R}$ of nontube type}
\label{table:irreducible_hermitian_symmetric_spaces}
\end{table}

On the level of Lie algebras, consider the Cartan decomposition $\mathfrak{g}^\textbf{R}=\mathfrak{h}^\textbf{R}+\mathfrak{m}^\textbf{R}$ and that of the maximal subtube $\mathfrak{g}_t^\textbf{R}=\mathfrak{h}_t^\textbf{R}+\mathfrak{m}_t^\textbf{R}$. There is an involution $\zeta:=\mathrm{Ad}(c^4):\mathfrak{g}^\textbf{R}\to \mathfrak{g}^\textbf{R}$ defined by a special element $c \in \mathfrak{h}^\textbf{R}+i\mathfrak{m}^\textbf{R}$ such that $\mathfrak{g}^\textbf{R}$ is of tube type if and only if $\zeta$ acts as the identity on $\mathfrak{g}^\textbf{R}$. If $\mathfrak{g}^\textbf{R}$ is of nontube type we can decompose the summands of the Cartan decomposition into the $\pm 1$-eigenspaces of $\zeta$ as 
\begin{align*}
\mathfrak{g}^\textbf{R}=\mathfrak{g}_T^\textbf{R}+\mathfrak{g}^\textbf{R}_-, \quad \mathfrak{h}^\textbf{R}=\mathfrak{h}_T^\textbf{R}+\mathfrak{h}^\textbf{R}_-, \quad \mathfrak{m}^\textbf{R}=\mathfrak{m}_t^\textbf{R}+\mathfrak{m}^\textbf{R}_-.
\end{align*}
where $\mathfrak{g}_T^\textbf{R}:=\mathfrak{n}_{\mathfrak{g}^\textbf{R}}(\mathfrak{g}_t^\textbf{R})$ and $\mathfrak{h}_T^\textbf{R}:=\mathfrak{n}_{\mathfrak{h}^\textbf{R}}(\mathfrak{h}_t^\textbf{R})$ are the normalizer of the maximal subtube and its maximal compact subalgebra.

On the group level, we denote by $G_T^\textbf{R}:=N_{G^\textbf{R}}(\mathfrak{g}_t^\textbf{R})_0$ and its maximal compact subgroup $H_T^\textbf{R}:=N_{H^\textbf{R}}(\mathfrak{h}_t^\textbf{R})_0$ the connected subgroups of $G^\textbf{R}$ with Lie algebras $\mathfrak{g}_T^\textbf{R}$ and $\mathfrak{h}_T^\textbf{R}$ normalizing the maximal subtube $G_t^\textbf{R}$ and $H_t^\textbf{R}$. Notice that since $\zeta$ acts as the identity on $\mathfrak{g}_T^\textbf{R}$, it is a Hermitian Lie algebra of tube type. In fact, $G_T^\textbf{R}$ is a maximal tube type Hermitian Lie subgroup of $G^\textbf{R}$, whose associated Hermitian symmetric space is irreducible if and only if $G_T^\textbf{R}=G_t^\textbf{R}$.

Consider the action of $H_T^\textbf{R}$ on $\mathfrak{m}_t^\textbf{R}$ by the adjoint action and let $\hat{C}^\textbf{R}:=H_T^\textbf{R} \cap \ker(\mathrm{Ad})|_{\mathfrak{m}_t^\textbf{R}}$ be the subgroup acting trivially on $\mathfrak{m}_t^\textbf{R}$. Then $G_T^\textbf{R}=G_t^\textbf{R}\times_{Z^\textbf{R}} \hat{C}^\textbf{R}$ is an extension of $G_t^\textbf{R}$ by $\hat{C}^\textbf{R}$ quotiented by the diagonal action of the central subgroup $Z^\textbf{R}=Z(G^\textbf{R})$ restricted to their common center. This fact will be useful for the proof of Theorem \ref{thm:maximal_components} and we will come back to it in Section \ref{sec:odd_magical_triples}. 

\begin{definition}
\label{def:toledo_character}
For a Hermitian Lie algebra $\mathfrak{g}^\textbf{R}$, let $\langle\cdot,\cdot \rangle$ be the Killing form on $\mathfrak{g}=\mathfrak{g}^\textbf{R}\otimes \textbf{C}$ and let $N$ be the dual Coxeter number of $\mathfrak{g}$. 
    The \textit{Toledo character} $\chi:\mathfrak{h}=\mathfrak{h}^\textbf{R}\otimes \textbf{C}\to \textbf{C}$ is defined by
\begin{align*}
    \chi(Y):=\frac{1}{N}\langle-iJ,Y\rangle,
\end{align*}
where $J\in \mathfrak{h}^\textbf{R}$ is the almost complex structure on $\mathfrak{m}^\textbf{R}=T_e(G^\textbf{R}/H^\textbf{R})$.
\end{definition}

In \cite[Proposition 2.5]{Biquard_García-Prada_Rubio_2017}, Biquard--Garc\'ia-Prada--Rubio give a criterion for when the Toledo character $\chi:\mathfrak{h}\to \textbf{C}$ lifts to a character $\tilde{\chi}:H\to \textbf{C}$. For our purposes it suffices to note that the Toledo character lifts to a character of $H$ for all of the simple Hermitian Lie groups of nontube type $G^\textbf{R}$ and their maximal subtubes $G^\textbf{R}_t$ in Table \ref{table:irreducible_hermitian_symmetric_spaces}. But in any case, a lift of the Toledo character is always possible up to an integer multiple of $\chi$ and we will call $\tilde{\chi}$ the \textit{lifted Toledo character}.

\subsubsection{Maximal $G^\textbf{R}$-Higgs bundles}
Let $G^\textbf{R}\subset G$ be a real form of a complex reductive Lie group with complexified Cartan decomposition $\mathfrak{g}=\mathfrak{h}+\mathfrak{m}$, and let $H\subset G$ be the complexification of the maximal compact subgroup $H^\textbf{R}\subset G^\textbf{R}$. Recall from our set up in Section \ref{sec:notations}, $X$ is a compact Riemann surface with genus $g\ge 2$ and we fixed a choice of $K^{1/2}$ for the canonical line bundle $K:=T^*X$. Let $L$ be a line bundle over $X$.

\begin{definition}
An \textit{$L$-twisted $G$-Higgs bundle} $\mathcal{E}_G$ is a pair $(E_G, \varphi)$, where $E_G\to X$ is a holomorphic $G$-principal bundle and $\varphi\in H^0(X, E_G[\mathfrak{g}]\otimes K)$ is a holomorphic section. Here $E_G[\mathfrak{g}]:=E_G\times_{\text{Ad}}\mathfrak{g}$ is the adjoint bundle.
\end{definition}

\begin{definition}
An \textit{$L$-twisted $G^\textbf{R}$-Higgs bundle} $\mathcal{E}_{G^\textbf{R}}$ is a pair $(E_H, \varphi)$, where $E_H\to X$ is a holomorphic $H$-principal bundle and $\varphi\in H^0(X, E_H[\mathfrak{m}]\otimes K)$ is a holomorphic section. Here $E_H[\mathfrak{m}]:=E_H\times_{\iota}\mathfrak{m}$ is the vector bundle associated to the isotropy representation $\iota:H\to \mathfrak{m}$ given by the adjoint action.
\end{definition}

The \textit{configuration space} of $L$-twisted $G$- and $G^\textbf{R}$-Higgs bundles, denoted by $\mathcal{H}_L(G)$ and $\mathcal{H}_L(G^\textbf{R})$, consists of the collection of Higgs pairs without imposing any stability conditions. Note that $\mathcal{H}_L(G^\textbf{R})\subseteq \mathcal{H}_L(G)$. \\

Given a $G$-Higgs bundle $\mathcal{E}_G$, the \textit{gauge group} $\mathcal{G}$ is the group of smooth automorphisms of the underlying bundle $E_G$ and acts on $\mathcal{H}(G)$ via pulling back the holomorphic structure and Higgs field. Similarly, given a $G^\textbf{R}$-Higgs bundle $\mathcal{E}_{G^\textbf{R}}$,the gauge group $\mathcal{G}_H$ of smooth automorphisms of the underlying bundle $E_H$ acts on $\mathcal{H}(G^\textbf{R})$ in an analogous manner. The orbits of the gauge group action are isomorphism classes of Higgs pairs.

We are ultimately interested in studying Higgs pairs in the configuration space up to isomorphism, or equivalently, up to the action of $\mathcal{G}$ on $\mathcal{H}_L(G)$ (resp. $\mathcal{G}_H$ on $\mathcal{H}_L(G^\textbf{R})$). However, crudely taking the quotient of the configuration space by the gauge group does not result in a nice Hausdorff space. It turns out that by throwing away certain \textit{unstable} objects and only considering the quotient of those that are \textit{polystable} does yield a tractable quotient called the \textit{moduli space of Higgs bundles} $\mathcal{M}(G)$ (resp. $\mathcal{M}(G^\textbf{R})$). To properly define the moduli space, we briefly recall the stability conditions of Higgs bundles.
\begin{definition}
    Let $\hat{G}\subset G$ be a complex reductive subgroup with Lie algebra $\hat{\mathfrak{g}}\subset \mathfrak{g}$. A $G$-Higgs bundle $\mathcal{E}_G=(E_G,\varphi)$ is said to reduce to a $\hat{G}$-Higgs bundle $\mathcal{E}_{\hat{G}}=(E_{\hat{G}},\varphi)$ if $E_G$ reduces to a holomorphic subbundle $E_{\hat{G}}$ such that $\varphi\in H^0(X,E_{\hat{G}}[\hat{\mathfrak{g}}]\otimes L)$.
\end{definition}

\begin{definition}
\label{def:real_reduction}
    Let $\hat{G}^\textbf{R}\subset G^\textbf{R}$ be a real reductive subgroup with maximal compact subgroup $\hat{H}^\textbf{R}\subset \hat{G}^\textbf{R}$ and complexified Cartan decomposition $\hat{\mathfrak{g}}=\hat{\mathfrak{h}}+\hat{\mathfrak{m}}$. A $G^\textbf{R}$-Higgs bundle $\mathcal{E}_{G^\textbf{R}}=(E_H,\varphi)$ is said to reduce to a $\hat{G}^\textbf{R}$-Higgs bundle $\mathcal{E}_{\hat{G}^\textbf{R}}=(E_{\hat{H}},\varphi)$ if $E_H$ reduces to a holomorphic subbundle $E_{\hat{H}}$ such that $\varphi \in H^0(X,E_{\hat{H}}[\hat{m}]\otimes L)$. 
\end{definition}

Let $P\subset G$ be a parabolic subgroup. A character $\chi_P:P\to \textbf{C}$ is called \textit{antidominant} if the associated line bundle $G(\chi_P):=G\times_{\chi_P}\textbf{C}\to G/P$ is ample. For a principal $P$-bundle $E_P$, we will denote by $E_P(\chi_P):=E_P\times_{\chi_P} \textbf{C}$ the line bundle associated to $\chi_P$.

\begin{definition}
A $G$-Higgs bundle $\mathcal{E}_{G}=(E_G,\varphi)\in \mathcal{H}_L(G)$ is:
\begin{itemize}
    \item \textit{stable} if for any parabolic subgroup $P\subset G$ and antidominant character $\chi_P$ and any reduction to $\mathcal{E}_P=(E_P,\varphi)$, we have that $\deg(E_{P}(\chi_P))> 0$,
    \item \textit{semistable} if for any parabolic subgroup $P\subset G$ and any antidominant character $\chi_P$ and any reduction to $\mathcal{E}_P=(E_P,\varphi)$, we have that $\deg(E_{P}(\chi_P))\ge 0$,
    \item \textit{polystable} if $\mathcal{E}_G$ is semistable and the reduction to $\mathcal{E}_P$ further reduces to an $L$-Higgs bundle $\mathcal{E}_L$, where $L$ is the Levi subgroup of $P$.
\end{itemize}
\end{definition}

The stability conditions of $G^\textbf{R}$-Higgs bundles are defined analogously for any parabolic subgroup $P^\textbf{R}$ of $G^\textbf{R}$ and the reduction to a $P^\textbf{R}$-Higgs bundle as per Definition \ref{def:real_reduction}.\\

Given the notion of stability, the \textit{moduli space of polystable $L$-twisted $G$- and $G^\textbf{R}$-Higgs bundles}
\begin{align*}
    \mathcal{M}_L(G):=\mathcal{H}^{ps}_L(G)/\mathcal{G}, \quad \mathcal{M}_L(G^\textbf{R}):=\mathcal{H}^{ps}_L(G^\textbf{R})/\mathcal{G}_H
\end{align*}
are defined as the space of gauge orbits of polystable Higgs bundles. This space can be obtained via a GIT construction, giving it the algebraic structure of a quasi-projective variety possibly with singularities \cite{Schmitt_2005}. Furthermore, by applying suitable Sobolev completions, the moduli space becomes a Hausdorff topological space \cite{Fan_2022}.  We will refer to this topology when discussing open and closed maps to the moduli space. However, since our results are built directly upon arguments in \cite{Bradlow_Collier_García-Prada_Gothen_Oliveira_2024}, we are spared from working with explicit descriptions of the topology of the moduli space.  \\

Now suppose $G^\textbf{R}$ is a simple Hermitian Lie group whose Toledo character $\chi$ lifts to a character $\tilde{\chi}$ on $H$. Given a $G^\textbf{R}$-Higgs bundle $\mathcal{E}_{G^\textbf{R}}:=(E_H,\varphi)$, let $E_H(\tilde{\chi}):=E_H\times_{\tilde{\chi}}\textbf{C}$ be the line bundle associated to the lifted Toledo character.
\begin{definition}
     The \textit{Toledo invariant} of $\mathcal{E}_{G^\textbf{R}}$ is defined to be the degree 
     \begin{align*}
         \tau(\mathcal{E}_{G^\textbf{R}}):=\deg (E_H(\tilde{\chi})).
     \end{align*}
    The Toledo invariant is an integer bounded above and below by the Milnor--Wood inequality \cite[Theorem 4.5]{Biquard_García-Prada_Rubio_2017}:
    $$|\tau| \le \mathrm{rank}(G^\textbf{R}/H^\textbf{R})(2g-2).$$ 
    A $G^\textbf{R}$-Higgs bundle $\mathcal{E}_{G^\textbf{R}}$ is called \textit{maximal} if its Toledo invariant  achieves maximal value,
$$\tau(\mathcal{E}_{G^\textbf{R}})=\mathrm{rank}(G^\textbf{R}/H^\textbf{R})(2g-2).$$

Furthermore, the Toledo invariant is a topological invariant of $\mathcal{M}(G^\textbf{R})$ that is constant on connected components. A connected component in $\mathcal{M}(G^\textbf{R})$ is called \textit{maximal} if it consists of maximal $G^\textbf{R}$-Higgs bundles. We denote the set of maximal $G^\textbf{R}$-Higgs bundles by $\mathcal{M}_{\mathrm{max}}(G^\textbf{R})\subset \mathcal{M}(G^\textbf{R})$.
\end{definition}

When $G^\textbf{R}$ is of nontube type, every maximal $G^\textbf{R}$-Higgs bundle is strictly polystable and the dimension of $\mathcal{M}_{\mathrm{max}}(G^\textbf{R})$ is strictly smaller than the expected dimension. We have the following statement:
\begin{theorem}[{\cite[Theorem 6.1, 6.2]{Biquard_García-Prada_Rubio_2017}}]
\label{thm:maximal_Higgs_bundles}
    Let $G^\textbf{R}$ be a simple Hermitian Lie group of nontube type. Recall that $G_T^\textbf{R}=G_t^\textbf{R}\times_{Z^\textbf{R}} {\hat{C}^\textbf{R}}$ is the maximal connected Hermitian Lie subgroup of $G^\textbf{R}$. Any maximal polystable $G^\textbf{R}$-Higgs bundle admits a reduction to a maximal polystable $G_T^\textbf{R}$-Higgs bundle. Furthermore, let $\mathrm{Ad}(G_t^\textbf{R})$ be the adjoint group of $G_t$, then there is a fibration of complex algebraic varieties
    \begin{align*}
        \mathcal{M}_0(\hat{C}^\textbf{R})\to \mathcal{M}_{\mathrm{max}}(G^\textbf{R})\to \mathcal{M}_{\mathrm{max}}(\mathrm{Ad}(G_t^\textbf{R}))
    \end{align*}
    where $\mathcal{M}_0(\hat{C}^\textbf{R})$ is the moduli space of $\hat{C}^\textbf{R}$-Higgs bundles with vanishing Toledo invariant.
\end{theorem}

The nonabelian Hodge correspondence \cite{Hitchin_1987,  Donaldson_1987, Simpson_1988, Corlette_1988, garciaprada2012hitchinkobayashicorrespondencehiggspairs}
establishes a homeomorphism between the moduli space of Higgs bundles to the character variety of surface group representations 
\begin{align*}
    \mathcal{M}(G^\textbf{R})\cong \mathcal{R}(G^\textbf{R})=\{\rho:\pi_1(X)\to G^\textbf{R}\}\sslash G^\textbf{R},
\end{align*}
taken as a GIT quotient with respect to the conjugation action of $G^\textbf{R}$.

Using bounded cohomology techniques, one can define the Toledo invariant for surface group representations into a Hermitian Lie group satisfying the same Milnor--Wood inequality \cite{Burger_Iozzi_Wienhard_2010}, which coincide with the definition for Higgs bundles under the nonabelian Hodge correspondence \cite[Proposition 4.2]{Biquard_García-Prada_Rubio_2017}, in other words, maximal $G^\textbf{R}$-Higgs bundles correspond to surface group representations $\pi_1(X)\to G^\textbf{R}$ with maximal Toledo invariant 
$$\mathcal{M}_{\mathrm{max}}(G^\textbf{R})\cong \mathcal{R}_{\mathrm{max}}(G^\textbf{R}).$$ 

Historically, Burger--Iozzi--Wienhard \cite{Burger_Iozzi_Wienhard_2010} first proved that any maximal representation $\pi_1(X)\to G^\textbf{R}$ factors through a maximal representation to $G_T^\textbf{R}$, while the Higgs bundle version in the form of Theorem \ref{thm:maximal_Higgs_bundles}, proven later with rather different techniques, is a manifestation of their result under the nonabelian Hodge correspondence. In particular, maximal representations of a nontube type Hermitian Lie group $G^\textbf{R}$ cannot be perturbed in a neighborhood of the character variety to have Zariski dense image. Such representations are called \textit{rigid}, i.e. contained in a neighborhood with no smooth homomorphisms, and this special phenomenon of maximal representations into a Hermitian Lie group of nontube type such that the maximal components to have dimension strictly less than the expected dimension $\dim \mathcal{M}_{\text{max}}(G^\textbf{R})< 2(g-1)\dim \mathfrak{g}^\textbf{R}$ is called \textit{rigidity}. On the other hand, a representation is called \textit{flexible} if it has a neighborhood in which smooth homomorphisms are dense.

\begin{remark}
\label{rmk:Kim-Pansu}
A result of Kim--Pansu \cite[Corollary 4]{Kim_Pansu_2014} shows that assuming the underlying curve $X$ has sufficiently large genus, surface group representations into a reductive real algebraic group are either flexible or rigid, and under this assumption, a representation is rigid only if it is a maximal representation onto its image, whose Zariski closure is a Hermitian Lie group of tube type. 
\end{remark}

\subsection{The Cayley correspondence}
Having described the rigidity phenomenon of maximal components when $G^\textbf{R}$ is Hermitian of nontube type, we turn to the description of maximal components when $G^\textbf{R}$ is Hermitian of tube type. The maximal components $\mathcal{M}_{\mathrm{max}}(G^\textbf{R})$ are captured by an isomorphism to the $K^2$-twisted moduli space of another group, which was first described systematically for classical Hermitian Lie groups in \cite{Bradlow_García-Prada_Gothen_2007} and was termed the \textit{Cayley correspondence} for its relation to the Cayley transformation of bounded symmetric domains. Later work by Bradlow--Collier--Garc\'ia-Prada--Gothen--Oliveira \cite{Bradlow_Collier_García-Prada_Gothen_Oliveira_2024} generalized the Cayley correspondence to real forms of complex Lie groups admitting a \textit{magical $\mathfrak{sl}_2$-triple}. 

We recount here the theory of magical $\mathfrak{sl}_2$-triples and that of the generalized Cayley correspondence with emphasis on its relation to surface group representations in higher Teichm\"uller theory.

\subsubsection{The Slodowy slice of an $\mathfrak{sl}_2$-triple}
We start by describing the uniformizing Higgs bundle over $X$. Let $S=\mathrm{SL}(2,\textbf{C})$, whose Lie algebra is generated by 
\begin{align*}
   \mathfrak{sl}_2\textbf{C}=\left\langle e=\begin{pmatrix} 0 & 1 \\ 0 & 0\end{pmatrix}, \,\, h=\begin{pmatrix} 1 & 0 \\ 0 & -1 \end{pmatrix}, \,\,f=\begin{pmatrix} 0 & 0 \\ 1 & 0\end{pmatrix} \right\rangle.
\end{align*}

Now consider $S^\textbf{R}=\mathrm{SL}(2,\textbf{R})$, whose complexified Cartan decomposition is given by $$\mathfrak{sl}_2\textbf{C}=\mathfrak{h}+\mathfrak{m}=\langle h \rangle+ \langle e,f\rangle.$$
Let $T\subset S$ be the complexified maximal compact subgroup obtained by the exponentiating the semisimple element $h$. The subgroup $T$ acts on $\mathfrak{m}$ by conjugation 
\begin{align*}
    \begin{pmatrix}
        \lambda & 0 \\ 0 & \lambda^{-1}
    \end{pmatrix} \cdot e=\lambda^2 e, \quad \begin{pmatrix}
        \lambda & 0 \\ 0 & \lambda^{-1}
    \end{pmatrix} \cdot f=\lambda^{-2} f.
\end{align*}

\begin{definition} 
Fixing a choice of $K^{1/2}$, the \textit{uniformizing $\mathrm{SL}(2,\textbf{R})$-Higgs bundle} is given by the $T$-Higgs bundle
$\mathcal{E}_T:=(E_T,\underline{f})$, where $E_T:=\accentset{\circ}{{K}}^{1/2}$ is the frame bundle of $K^{1/2}$ and $\underline{f}\in H^0(X,E_T[\langle e,f\rangle]\otimes K)$ is the constant section valued at 1 in 
$$ H^0(X,E_T[\langle f\rangle]\otimes K)\cong H^0(X,\mathcal{O}_X),$$ 
where the last isomorphism uses the fact that $E_T[\langle f \rangle] \cong K^{-1}$ that can be deduced by the action of $T$ on $f$. The vector bundle associated to $\mathcal{E}_T$ via the standard representation takes the familiar form
\begin{align*}
\mathcal{E}_T[\textbf{C}^2]=\left(E_T[\textbf{C}^2]=K^{1/2}\oplus K^{-1/2},\underline{f}=\begin{pmatrix}
        0 & 0 \\ 1 & 0
    \end{pmatrix}\in H^0(X,E_T[\textbf{C}^2]\otimes K) \right)
\end{align*}
\end{definition}

Let $\rho:\mathfrak{sl}_2\textbf{C}\to \mathfrak{s}=\langle e,h,f \rangle \subset \mathfrak{g}$ be an $\mathfrak{sl}_2$-triple of $G$ and let $S\subset G$ be connected Lie subgroup associated to the image of $\mathfrak{sl}_2$ with maximal torus $T\subset S \subset G$. Let $C \subset G$ be the centralizer of $\mathfrak{s}$ with Lie algebra $\mathfrak{c}$ and recall that $V_j \subset W_j$ are the highest weight spaces of the $\mathfrak{sl}_2$-module decomposition $\mathfrak{g}=\bigoplus_j W_j$ determined by $\rho$. The centralizer $C$ stabilizes the sum of nontrivial highest weight spaces $V_\rho^{>0}=\bigoplus_{j>0} V_j$.

\begin{definition}
    The \textit{$\rho$-Slodowy category} $\mathcal{B}_\rho(G)$ consists of the following data:
    \begin{itemize}
        \item \textit{Objects} are tuples $(\mathcal{E}_C,\{\varphi_{j}\})$ where $\mathcal{E}_C=(E_C,\varphi)$ is a $C$-Higgs bundle and $\varphi_{j}\in H^0(X, \mathcal{E}_C[V_j]\otimes K)$ are holomorphic sections taking values in the nontrivial highest weight spaces $V_{j>0}$.
        \item \textit{Morphisms} between two objects $(\mathcal{E}_C,\{\varphi_j\})$ and $(\mathcal{E}'_C,\{\varphi_j'\})$ are isomorphisms $\Phi:\mathcal{E}_C\xrightarrow[]{\cong} \mathcal{E}_C'$ of $C$-Higgs bundles such that $\varphi_j=\Phi^*\varphi_j'$ for all $j$.
    \end{itemize}
\end{definition}

\begin{definition}
The \textit{Slodowy map} sends an objects in the Slodowy category to a $G$-Higgs pair 
\begin{align}
\label{eqn:slodowy_map}
\hat{\Psi}_\rho:\,\,&\mathcal{B}_\rho(G) \longrightarrow \mathcal{H}(G), \\
    &((E_C,\varphi_C),\{\varphi_j\}) \mapsto (E_G:=(E_C\times E_T)[G], \,\,f+\varphi_C+\sum_{j} \varphi_{j}),
\end{align}
where $(E_T,\underline{f})$ is the uniformizing Higgs bundle and $E_G$ is extension of the $C\times T$-bundle to a $G$-bundle. Moreover, since $C$ acts trivially on $f$, the section $\underline{f} \in H^0(X,E_T[\langle f\rangle ]\otimes K)$ may be considered as a holomorphic section of $H^0(X,E_G[\mathfrak{g}]\otimes K)$. 
\end{definition}

\begin{definition}
Given a real form $G^\textbf{R}\subset G$, consider the restriction of the Slodowy map $\hat{\Psi}_\rho$ to the preimage of polystable $G^\textbf{R}$-Higgs pairs in $\mathcal{H}(G^\textbf{R})$. Upon descending to the moduli space by identifying gauge orbits, we will call its image in $\mathcal{M}(G^\textbf{R})$ the \textit{Slodowy slice} $\mathrm{Slo}_\rho$ associated to $\rho$.
\end{definition}

\begin{remark}
The Slodowy map we describe here is a simplification of the Slodowy functor for a fixed uniformizing Higgs bundle $\Theta=(E_T,\underline{f})$, where we forget the $(G,P)$-oper structure,
\begin{align*}
\hat{\Psi}_{\rho,\Theta}:\,\,&\mathcal{B}_\rho(G) \longrightarrow \mathcal{O}p^0(G,P_\rho), \\
    &((E_C,\varphi_C),\{\varphi_j\}) \mapsto (E_G,\, E_{P_\rho},\, \underline{f}+\varphi_C+\sum_{j=1}^N \varphi_{j}).
\end{align*}
In \cite{Collier_Sanders_2021}, Collier--Sanders showed that this map is an equivalence of categories between the Slodowy category and $(G,P_\rho)$-opers when $\lambda=0$.
\end{remark}

\subsubsection{Magical $\mathfrak{sl}_2$-triples and the Cayley correspondence}
Let $\mathfrak{g}$ be a simple complex Lie algebra and $\mathfrak{g}^\textbf{R}\subset \mathfrak{g}$ a real form.
\begin{definition}
    An $\mathfrak{sl}_2$-triple $\rho:\mathfrak{sl}_2\to \mathfrak{g}$ is called \textit{magical} if the \textit{magical involution} $\sigma_\rho:\mathfrak{g}\to \mathfrak{g}$
\begin{align}
\label{eqn:magical_involution}
\sigma_\rho(x)=\begin{cases}x, & x\in V_0, \\ (-1)^{k+1}x, & x\in (\mathrm{ad}_f)^k\,(V_j), \,\, j>0,\end{cases}
\end{align}
    which is apriori only a vector space involution, is a Lie algebra involution. 
\end{definition}

Given a magical triple $\rho$, there is a real form $\tau_\rho$  associated to $\sigma_\rho$ such that $\tau_\rho \sigma_\rho=\sigma_\rho \tau_\rho$, whose fixed points $
\mathfrak{g}^{\tau_\rho}=\mathfrak{g}_\rho^\textbf{R}$ defines a real form of $\mathfrak{g}$ which we call the \textit{canonical real form of $\rho$}.
\begin{definition}
    A real $\mathfrak{sl}_2$-triple $\hat{\rho}:\mathfrak{sl}_2\textbf{R}\to \mathfrak{g}^\textbf{R}$ is called \textit{magical} if its Cayley transform $\rho:\mathfrak{sl}_2\textbf{C}\to \mathfrak{g}$ is magical and has $\mathfrak{g}^\textbf{R}=\mathfrak{g}^\textbf{R}_\rho$ as its canonical real form. 
\end{definition}

\begin{example}
\label{principal}
The principal $\mathfrak{sl}_2$-triple $\rho_{\text{prin}}:\mathfrak{sl}_2\textbf{C}\to \mathfrak{sl}_3\textbf{C}$ is a magical triple:
\begin{align*}
        \rho_{\text{prin}}=\left\langle e=\begin{pmatrix}
            0 & 1 & 0\\ 
            0 & 0 & 1 \\
            0 & 0 & 0
        \end{pmatrix}, \quad 
        h=\begin{pmatrix}
            2 & 0 & 0 \\ 
            0 & 0 & 0 \\
            0 & 0 & -2
        \end{pmatrix}, \quad
        f=\begin{pmatrix}
            0 & 0 & 0\\ 
            2 & 0 & 0 \\
            0 & 2 & 0
        \end{pmatrix}\right\rangle.
    \end{align*}

The $\mathfrak{sl}_2$-module decomposition of $\mathfrak{sl}_2\textbf{C}$ with respect to $\rho$ can be presented diagrammatically as follows:
\begin{center}
\begin{tabular}{ c c c c c c}
$W_4$ & $\bullet\spc^{\bigminus}$ & $\bullet\spc^{\bigplus}$ & $\bullet\spc^{\bigminus}$ & $\bullet\spc^{\bigplus}$ & $\bullet\spc^{\bigminus}$ \\ 
 $W_2$ &  & $\bullet\spc^{\bigminus}$ & $\bullet\spc^{\bigplus}$ & $\bullet\spc^{\bigminus}$ & \\  
  & $\mathfrak{g}_{-4}$ & $\mathfrak{g}_{-2}$ & $\mathfrak{g}_{0}$ & $\mathfrak{g}_{2}$ & $\mathfrak{g}_{4}$   
\end{tabular}
\end{center}
where each dot represents the intersection $W_j\cap \mathfrak{g}_{j-2k}$ of the two decompositions $\mathfrak{g}=W_2\oplus W_4$ and $\mathfrak{g}=\mathfrak{g}_{-4}\oplus \mathfrak{g}_{-2}\oplus \mathfrak{g}_{0}\oplus \mathfrak{g}_{2}\oplus \mathfrak{g}_{4}$ described in Section \ref{sec:sl2}. The action of $\sigma_\rho$ on $\mathfrak{g}$ is illustrated by its $\pm 1$-eigenspaces as labeled. The condition for $\rho_{\text{prin}}$ to be magical is for $\sigma_\rho$ to be a Lie algebra involution, which is satisfied in this case. The canonical real form of $\rho_{\text{prin}}$ is the split real form $\mathfrak{g}^\textbf{R}_\rho=\mathfrak{sl}_3\textbf{R}$.
\end{example}

Magical $\mathfrak{sl}_2$-triples $\hat{\rho}:\mathfrak{sl}_2\textbf{R}\to \mathfrak{g}^\textbf{R}$ are classified in \cite[Theorem 3.1]{Bradlow_Collier_García-Prada_Gothen_Oliveira_2024} by the weighted Dynkin diagram of their Cayley transforms $\rho:\mathfrak{sl}_2\textbf{C} \to \mathfrak{g}$. The following is the list of real forms $\mathfrak{g}^\textbf{R}$ admitting magical triples:
\begin{theorem}[\cite{Bradlow_Collier_García-Prada_Gothen_Oliveira_2024}, Theorem C]
    \label{thm:classification} 
Let $\mathfrak{g}$ be a simple complex Lie algebra with real form $\mathfrak{g}^\textbf{R}\subset \mathfrak{g}$. Then $\mathfrak{g}^\textbf{R}$ is canonical real forms of some magical $\mathfrak{sl}_2$-triple if and only if it is in one of the four families:
\begin{enumerate}
    \item [(1)] $\mathfrak{g}$ is split real and $\mathfrak{g}^\textbf{R}$ is the split real form,
    \item [(2)] $\mathfrak{g}$ is of type $A_{2n-1}, B_n, C_n, D_n, D_{2n}$ or $E_7$ and $\mathfrak{g}^\textbf{R}$ is Hermitian of tube type,
    \item [(3)] $\mathfrak{g}=\mathfrak{so}_N\textbf{C}$ and $\mathfrak{g}^\textbf{R}=\mathfrak{so}_{p,N-p}$, $p\ge 3$,
    \item [(4)] $\mathfrak{g}$ is exceptional and $\mathfrak{g}^\textbf{R}$ is one of the following:
    $\begin{array}{c| c| c| c| c}
        \mathfrak{g} & E_6 & E_7& E_8 &F_4  \\
        \hline
        \mathfrak{g}^\textbf{R} & E_6^2 & E_7^{-5}& E_8^{-24}& F_4^4
    \end{array}.$
\end{enumerate}
\end{theorem}

\begin{remark}
\label{rmk:even_magical_triple}
Notably, all of the magical triples are even, i.e. the eigenvalues of $\mathrm{ad}_h$ are all even, or equivalently, the $\mathrm{\textbf{Z}}$-grading $\mathfrak{g}=\bigoplus \mathfrak{g}_{2j}$ consists entirely of even weight spaces.
\end{remark}

\begin{example}
We list here the magical triples of the Hermitian Lie algebras appearing as the maximal subtube of a Hermitian Lie algebra of nontube type in Table \ref{table:irreducible_hermitian_symmetric_spaces}. The label on each root is the weight $\mathrm{ad}_h(\alpha)$, from which their associated $\mathfrak{sl}_2$-data can be calculated \cite[Section 4]{Bradlow_Collier_García-Prada_Gothen_Oliveira_2024}.

\begin{align*}
&\mathfrak{g}^\textbf{R}=\mathfrak{su}(p,p), \,\,\mathfrak{g}=\mathfrak{sl}_{2p}\textbf{C}, \,\,\,
    &&A_{2p-1}: \,\,\dynkin[mark=o, make indefinite edge/.list={1-2,4-5}, labels={0,0,2,0,0}, labels*={,,\alpha_{p},,}]A5 \\
&\mathfrak{g}^\textbf{R}=\mathfrak{so}^*_{4m}, \,\,\mathfrak{g}=\mathfrak{so}_{4m}\textbf{C}, \,\,\,
    &&D_{2(2n+1)}: \,\,\dynkin[mark=o, labels={0,0,0,0,2,0}]D{}\\
   & &&D_{2(2n)}: \,\,\dynkin[mark=o, labels={0,0,0,0,0,2}]D{}\\
    &\mathfrak{g}^\textbf{R}=\mathfrak{so}(2,8), \,\,\mathfrak{g}=\mathfrak{so}_{10}\textbf{C}, \,\,\,
    &&D_5:\,\, \dynkin[mark=o, labels={2,0,0,0,0}]D5
    \end{align*}
\end{example}    

\begin{remark}
\label{rmk:relation_to_theta_positivity}
    The theory of magical $\mathfrak{sl}_2$-triples are closely related to that of $\Theta$-positive structures defined in \cite{Guichard_Wienhard_2025}. In particular, real form $G^\textbf{R} \subset G$ admits a magical triple if and only if it admits a $\Theta$-positive structure \cite[Theorem 8.14]{Bradlow_Collier_García-Prada_Gothen_Oliveira_2024}. Concretely, if $G^\textbf{R}$ admits a $\Theta$-positive structure, the $\mathfrak{sl}_2$-triple induced by the $\Theta$-principal embedding in \cite[Section 7.1]{Guichard_Wienhard_2025} is magical. Conversely, given a magical $\mathfrak{sl}$-triple $\hat{
    \rho
    }$, consider the set of weight-2 simple roots in weighted Dynkin diagram of its Cayley transform $\rho$. Folding the Dynkin diagram with respect to the involution defining $\mathfrak{g}^\textbf{R}$, we obtain the Dynkin diagram of the restricted root system. The simple restricted roots corresponding to the weight-2 roots then forms the subset $\Theta\subset \Delta$ giving rise to a $\Theta$-positive structure. 
\end{remark}

\begin{definition}
For a magical triple $\rho:\mathfrak{sl}_2\to \mathfrak{g}$, there is a Lie algebra involution $\theta_\rho:\mathfrak{g}_0\to \mathfrak{g}_0$ given by
\begin{align}
\label{eq:theta_involution}
    \theta_\rho(x)=
    \begin{cases}
    x, \quad \quad &x\in \mathfrak{c} \\
    -x, \quad & \text{otherwise}.
    \end{cases}
\end{align}
The real form $\tau_0$ associated to $\theta_\rho$ fixes a real form $\mathfrak{g}_\mathcal{C}^\textbf{R}$ of $\mathfrak{g}_0$ called the \textit{Cayley real form}. 
\end{definition}

Now to each magical triple $\rho$, there exists at most one weight $m_c$ such that the dimension of the highest weight space $\dim V_{m_c}\ge 1$ \cite[Proposition 4.3]{Bradlow_Collier_García-Prada_Gothen_Oliveira_2024}. Let $\{l_j\}$ be the weights of the rest of the highest weight spaces such that $\mathfrak{g}=W_{m_c}\bigoplus_{j=1}^{r(\rho)} W_{l_j}$. With this set up, the generalized Cayley correspondence for a magical triple $\hat{\rho}$ is the restriction of the Slodowy map $\hat{\Psi}_\rho$ descended to the level of moduli spaces:

\begin{theorem}[\cite{Bradlow_Collier_García-Prada_Gothen_Oliveira_2024}, Theorem 7.1]
\label{thm:cayley_correspondence}
    Let $G$ be a complex Lie group and let $\hat{\rho}$ be a magical triple of the real form $G^\textbf{R}\subset G$. Let  $\tilde{\mathfrak{g}}^\textbf{R}$ to be the semisimple part of $\mathfrak{g}_\mathcal{C}^\mathrm{\textbf{R}}$, whose complexified maximal compact subalgebra is $\mathfrak{c}\subset \mathfrak{g}_0$. On the group level, we define $\tilde{G}^\textbf{R}$ to be the subgroup of $G^\textbf{R}$ with Lie algebra $\tilde{\mathfrak{g}}^\textbf{R}$ and maximal compact subgroup $\mathcal{C}:=C\cap G^\textbf{R}$. Then the \textit{Cayley map} given by
    \begin{align*}
        \Psi_\rho:\,\, &\mathcal{M}_{K^{m_c+1}}(\tilde{G}^\mathrm{\textbf{R}})\times\bigoplus_{j=1}^{r(\rho)}H^0(K^{l_j+1}) \,\,\longrightarrow \,\,\mathcal{M}(G^\mathrm{\textbf{R}}), \\ &((E_\mathcal{C},\psi_{m_c}),\varphi_1,...,\varphi_{r(\rho)}) \mapsto (E_H:=(E_\mathcal{C}\times_X E_T)[H], \,\,\underline{f}+\varphi_{m_c}+\sum_{j=1}^{r(\rho)} \varphi_j),
    \end{align*}
    where $\varphi_{m_c}:=\mathrm{ad}_f^{-m_c}(\psi_{m_c}) \in H^0(X,E_H[V_{m_c}]\otimes K)$, is well-defined on the moduli spaces that is injective, open and closed.
\end{theorem}

\begin{example}
\label{ex:Cayley_correspondence_Hermitian}
For Hermitian Lie groups $G^\textbf{R}$ of tube type, we have $m_c=2$ and $\{l_j\}=\{0\}$, and the Cayley correspondence is an isomorphism onto the maximal components
\begin{align*}
    \Psi_\rho:\,\, &\mathcal{M}_{K^2}(\tilde{G}^\mathrm{\textbf{R}})\times H^0(K) \xlongrightarrow{\cong} \mathcal{M}_{\mathrm{max}}(G^\textbf{R})\subset \mathcal{M}(G^\mathrm{\textbf{R}}) 
\end{align*}
If we denote $H^*\subset H$ the noncompact dual of $H^\textbf{R}$, we can subsume the section $H^0(K)$ and the Cayley correspondence becomes the isomorphism 
\begin{align*}
 \mathcal{M}_{K^2}(H^*)\cong \mathcal{M}_{\mathrm{max}}(G^\textbf{R}).
\end{align*}
This phenomenon is the precursor to the generalized Cayley correspondence.
\end{example}

\begin{remark}
In recent joint work with E. Chen and M. Yang \cite{chen2025baabraneshigherteichmullertheory}, we extended the generalized Cayley correspondence to an open and universally closed Cayley morphism of $BAA$-branes on the level of Higgs stacks by employing the notion of Gaiotto's Lagrangian on two suitable Hamiltonian $G$-spaces.
\end{remark}

\section{Extended magical triples and their classification}
Let $\mathfrak{g}$ be a simple complex Lie algebra and $\mathfrak{g}^\textbf{R}$ a real form.

\begin{definition}
\label{def:extended_magical_triple}
A complex $\mathfrak{sl}_2$-triple $\rho:\mathfrak{sl}_2\textbf{C}\to \mathfrak{g}$ is called an \textit{extended magical triple} if there exists an involution $\sigma_\rho:\mathfrak{g}\to \mathfrak{g}$ that agrees with the magical involution \eqref{eqn:magical_involution} on the even weight spaces of the $\mathfrak{sl}_2\textbf{C}$-module decomposition $\mathfrak{g}=W_{\text{odd}}\oplus W_{\text{even}}$. If $\rho$ is magical, the real form $\mathfrak{g}_\rho^\textbf{R}$ associated to $\sigma_\rho$ will be called the \textit{canonical real form}. 
\end{definition}

\begin{definition}
A real $\mathfrak{sl}_2$-triple $\hat{\rho}:\mathfrak{sl}_2\textbf{R}\to \mathfrak{g}^\textbf{R}$ is said to be \textit{extended magical} if its Cayley transform $\rho$ is extended magical with canonical real form $\mathfrak{g}_\rho^\textbf{R}=\mathfrak{g}^\textbf{R}$. 
\end{definition}

Let $\sigma: \mathfrak{g}\to \mathfrak{g}$  be Lie algebra involution giving rise to the complexified Cartan decomposition $\mathfrak{g}=\mathfrak{h}+\mathfrak{m}$. 
We can equivalently define a real $\mathfrak{sl}_2$-triple $\hat{\rho}: \mathfrak{sl}_2\textbf{R}\to \mathfrak{g}^\textbf{R}$ to be an extended magical triple if, with respect to the $\mathfrak{sl}_2\mathbf{C}$-module decomposition of its Cayley transform $\rho:\mathfrak{sl}_2\mathbf{C}\to \mathfrak{g}$, $\sigma$ acts as $+\mathrm{Id}$ on $\mathfrak{c}$ and as $-\mathrm{Id}$ on the nontrivial even highest weight spaces $V_{even}:=\bigoplus_{j >0} V_{2j}$. In this case, $\sigma$ will be the involution for $\rho$ that agrees with the magical involution on even weight spaces. 

Recall from our set up that $V_j=\bigoplus_j {v_j}^{n_j}$ is the sum of the highest weight spaces $v_j$ from the $n_j$-copies of $w_j$. Let each $v_j$ be an $\pm 1$-eigenspace of $\sigma$. Since $\rho$ is a normal triple, i.e. $\sigma(f)=-f$, the sign of $\sigma$ on $(\mathrm{ad}_f)^k (v_j)$ alternates as we apply the lowering operator $\mathrm{ad}_f$, and the involution $\sigma:\mathfrak{g}\to \mathfrak{g}$ is completely determined by its value on $v_j$. By Definition \ref{def:extended_magical_triple}, $\rho$ is an extended magical triple if and only if
     \begin{align}
        \sigma(x)=
        \begin{cases}
            x, & x \in \mathfrak{c},\\
            (-1)^{k+1}x, & x \in (\mathrm{ad}_f)^k(V_{2j>0}).
        \end{cases}
    \end{align}
    
One observes immediately that the definition of an extended magical triple agrees with that of a magical triple on the set of even $\mathfrak{sl}_2$-triples. Therefore, by Remark \ref{rmk:even_magical_triple}, an even $\mathfrak{sl}_2$-triple is an extended magical triple if and only if it is a magical triple. Any extended magical triple that is not on the list of magical triples \cite[Theorem 3.1]{Bradlow_Collier_García-Prada_Gothen_Oliveira_2024} must therefore be odd.

\begin{definition}
An extended magical triple that is magical in the sense of \cite{Bradlow_Collier_García-Prada_Gothen_Oliveira_2024} will be called an \textit{even magical triple}. Otherwise it will be called an \textit{odd magical triple}. 
\end{definition} 

Much as in the case for magical triples \cite[Proposition 2.11]{Bradlow_Collier_García-Prada_Gothen_Oliveira_2024}, there is a straightforward criterion for detecting extended magical triples:
\begin{proposition}
\label{prop:modified_criteria}
An $\mathfrak{sl}_2$-triple $\hat{\rho}:\mathfrak{sl}_2\mathrm{\textbf{R}}\to \mathfrak{g}^\textbf{R}$ is an extended magical triple if and only if $\mathfrak{c}\subset \mathfrak{h}$ and 
\begin{align}
\label{eqn:extended_triple_criteria}
    \dim \mathfrak{m}-\dim \mathfrak{h} =\dim \mathfrak{g}_0-2\dim \mathfrak{c}.
\end{align}
In particular, if $\hat{\rho}$ is an extended magical triple, then the centralizer $\mathfrak{c}^\textbf{R}\subset \mathfrak{h}^\textbf{R}$ is a subalgebra of the maximal compact subalgebra and is compact, i.e. the Killing form restricted to $\mathfrak{c}^\textbf{R}$ is negative definite.
\end{proposition}

\begin{proof}
To begin, recall that each $v_j$ is a $\pm1$-eigenspace of $\sigma$, and $\sigma:\mathfrak{g}\to \mathfrak{g}$ is completely determined by its values on $v_j$ in the sense that $$\mathrm{sgn}(\sigma((\mathrm{ad}_f)^k\,v_j))=(-1)^k\mathrm{sgn}(\sigma(v_j)).$$

We claim that the odd weight spaces $w_{2j+1}$ do not contribute to the difference
\begin{align*}
    \dim \mathfrak{m}-\dim \mathfrak{h}&=\sum_{j=0}^M \dim (\mathfrak{m}\cap W_{j})-\dim (\mathfrak{h}\cap W_{j})\\
    &= \sum_{j=0}^M \dim (\mathfrak{m}\cap V_{2j})-\dim (\mathfrak{h}\cap V_{2j}).
\end{align*}

This can be observed by computing the difference in the number of $+1$ and $-1$-eigenspaces
\begin{align*}
    \dim (\mathfrak{m}\cap w_j)-\dim (\mathfrak{h}\cap w_j)=\begin{cases}-\mathrm{sgn}(\sigma(v_j)), & j \,\, \text{even},\\
    0, & j \,\, \text{odd}.\end{cases}
\end{align*}

Suppose $\hat{\rho}$ is an extended magical triple, then $\mathfrak{c}$ is fixed by $\sigma$ and therefore is contained in $\mathfrak{h}$. Moreover, the nontrivial even highest weight spaces $V_{even}:=\bigoplus_{j>0} V_{2j}$ are $-1$-eigenspaces of $\sigma$. Therefore, $$\dim \mathfrak{m}-\dim \mathfrak{h}= \dim \bigoplus_{j>0} V_{2j}-\dim \mathfrak{c}=\dim \mathfrak{g}_0-2\dim \mathfrak{c},$$
where we note that the dimension of the nontrivial even highest weight spaces is $\dim \mathfrak{g}_0-\dim \mathfrak{c}$.

For the opposite direction, suppose $\mathfrak{c} \subset \mathfrak{h}$ and \eqref{eqn:extended_triple_criteria} holds. Then
\begin{align*}
    \dim \mathfrak{m}-\dim \mathfrak{h}&= \sum_{j=0}^M \dim (\mathfrak{m}\cap V_{2j})-\dim (\mathfrak{h}\cap V_{2j})\\
    &=\sum_{j=1}^M \dim (\mathfrak{m}\cap V_{2j})-\dim (\mathfrak{h}\cap V_{2j})-\dim \mathfrak{c}\\
    &=\dim \mathfrak{g}_0-2\dim \mathfrak{c}.
\end{align*}
Hence, we have that
\begin{align*}
    \sum_{j=1}^M \dim (\mathfrak{m}\cap V_{2j})-\dim (\mathfrak{h}\cap V_{2j})&=\dim \mathfrak{g}_0-\dim \mathfrak{c} = \sum_{j=1}^M \dim V_{2j}.
\end{align*}
This is only possible when the left hand side attains its maximum, i.e. if and only if the nontrivial even highest weight spaces $V_{even} \subset \mathfrak{m}$. This concludes the proof.
\end{proof}

\begin{theorem}
\label{thm:classification_theorem}
Let $\mathfrak{g}^\textbf{R}$ be the real form of a simple complex Lie algebra $\mathfrak{g}$. Let $\hat{\rho}:\mathfrak{sl}_2\mathrm{\textbf{R}} \to \mathfrak{g}^\textbf{R}$ be an $\mathfrak{sl}_2$-triple and $\rho:\mathfrak{sl}_2\mathbf{C} \to \mathfrak{g}$ its Cayley transform. Then $\hat{\rho}$ is an odd magical triple if and only if $\mathfrak{g}^\textbf{R}$ is a Hermitian Lie algebra of nontube type and the weighted Dynkin diagram of $\rho$ is one of the following: 
\begin{align*}
&\mathfrak{g}^\textbf{R}=\mathfrak{su}(p,q), \,\,\mathfrak{g}=\mathfrak{sl}_{p+q}\mathrm{\textbf{C}}, \,\,\,
    &&A_{p+q-1}: \,\,\dynkin[mark=o, make indefinite edge/.list={1-2,4-5,7-8}, labels={0,0,1,0,0,1,0,0}, labels*={,,\alpha_{p},,,\alpha_{q}}]A8 \\
&\mathfrak{g}^\textbf{R}=\mathfrak{so}^*_{4m+2}, \,\,\mathfrak{g}=\mathfrak{so}_{4m+2}\mathrm{\textbf{C}}, \,\,\,
    &&D_{2n+1}: \,\,\dynkin[mark=o, labels={0,0,0,0,1,1}]D{}\\
    &\mathfrak{g}^\textbf{R}=E_6^{-14}, \,\,\mathfrak{g}=E_6, \,\,\,
    &&E_6:\,\, \dynkin[mark=o, labels={1,0,0,0,1,0}]E6
\end{align*}
\end{theorem}

The proof relies on the following lemma:
\begin{lemma}
\label{lem:dimension_of_g0}
    Let $\mathfrak{g}$ be a classical Lie algebra and $\rho:\mathfrak{sl}_2\mathbf{C}\hookrightarrow \langle e,h,f\rangle \subset \mathfrak{g}$ an $\mathfrak{sl}_2$-triple associated to the partition $\sum_i r_i\cdot i$. Then the dimension of $\mathfrak{g}_0$ is given by 
    \begin{align*}
        \dim \mathfrak{g}_0 =\dim V_\rho- \begin{cases}
        \sum_{\substack{i< j \\
        i\neq j \,(\mathrm{mod} \,2)}}\, 2i\,r_ir_j, & \mathfrak{g}=\mathfrak{sl}_n\mathrm{\textbf{C}} \\
        \sum_{\substack{i< j \\
        i\neq j \,(\mathrm{mod} \,2)}}\, i\, r_ir_j, & \mathfrak{g}=\mathfrak{sp}_{2n}\mathrm{\textbf{C}} \\
        \sum_{\substack{i< j \\
        i\neq j \,(\mathrm{mod} \,2)}}\, i\, r_ir_j, & \mathfrak{g}=\mathfrak{so}_N\mathrm{\textbf{C}}
        \end{cases}
    \end{align*}
    In particular, $\rho$ is an even triple, i.e. $\dim \mathfrak{g}_0=\dim V_\rho$ if and only if $r_i=0$ for all even $i$ or $r_i=0$ for all odd $i$.
\end{lemma}

\begin{proof}
The proof of Theorem \ref{thm:characterization} can be abstracted from the proof of Theorem 6.1.3 in \cite{Collingwood_McGovern_1993}. We follow their strategy closely and adopt their conventions. Let $V$ be the standard representation of $\mathfrak{g}$. We have an $\mathfrak{sl}_2\mathbf{C}$-module decomposition
\begin{align}
\label{eqn:standard_rep_decomposition}
    V\cong \bigoplus_{d\ge 0} M(d)^{\oplus \,r_{d+1}}
\end{align}
where the number of irreducible summands $M(d)$  corresponds to the partition $\textbf{r}=[r_1,...,r_N]$ associated to $\rho$. Let $L(d)$ denote the lowest weight space of $M(d)$ and let $S_{d,e}$ denote the space of linear maps $L(d)\to M(e)$ sending the lowest weight vector $v$ in $L(d)$ to a vector $w$ in $M(e)$ such that $(\mathrm{ad}_e)^{d+1}(w)=0$. Theorem 6.1.3 asserts that that sum of the highest weight spaces of $\mathfrak{g}$, $$V_\rho=\mathrm{ker}(\mathrm{ad}_e)=\{X\in \mathfrak{g}, \mathrm{ad}_e(X)=X\}$$ is isomorphic to the sum of the projections $S_{d,e}$, up to certain counting factors occuring from extra symmetries when $\mathfrak{g}$ is $\mathfrak{so}_N\mathrm{\textbf{C}}$ or $\mathfrak{sp}_{2n}\mathrm{\textbf{C}}$. From this, the authors were able to derive a formula for the dimension of $V_\rho$, which we include here for completeness:
\begin{align*}
    \dim V_\rho=\begin{cases} \sum_i s_i^2 -1, & \mathfrak{g}=\mathfrak{
    sl
    }_n\textbf{C} \\
    \frac{1}{2}\sum_i s_i^2 +\frac{1}{2}\sum_{i\, \text{odd}} r_i^2, & \mathfrak{g}=\mathfrak{sp}_{2n}\textbf{C}\\
    \frac{1}{2}\sum_i s_i^2 -\frac{1}{2}\sum_{i\, \text{odd}} r_i^2, & \mathfrak{g}=\mathfrak{so}_{N}\textbf{C}\\
    \end{cases}
\end{align*} 
where $\sum_i s_i=n$ is the dual partition of $\sum_i i\cdot r_i$.

For our purposes, note that $\mathfrak{g}_0$ is isomorphic to the even highest weight spaces by applying the raising operator $\mathrm{ad}_e$. We show the following a linear map in $S_{d,e}$ corresponds to a highest weight vector in an odd irreducible representation $W_{2j+1}$ of $\mathfrak{g}$ if and only if $d\neq e\,\, (\mathrm{mod} \,\,2)$. Consequently, the dimension of $\mathfrak{g}_0$ differs from $\dim V_\rho$ precisely where $d$ and $e$ of opposite parities contributes. 

Let $X$ be a highest weight vector in $S_{d,e}$. Following a description of the semisimple element $h$ for the various classical Lie algebras in Section 5.2 of \cite{Collingwood_McGovern_1993}, we note that with respect to the decomposition \eqref{eqn:standard_rep_decomposition}, $h$ takes the following form:
\begin{align*}
    h=\begin{pmatrix}D(r_1)& 0 & \cdots & 0 \\
    0 & D(r_2) & \cdots & 0 \\
    \vdots & 0 & \ddots & 0 \\
    0 & \cdots & 0 & D(r_N)\end{pmatrix}
\end{align*}
where each $D(r_i)$ is a diagonal matrix with integer entries that are of opposite parity to that of $r_i$, i.e. if $r_i$ is even, then $D(r_i)$ has odd entries and vice versa. The weight of $X$ under $\mathrm{ad}_h$ is the difference of two respective entries in $D(d)$ and $D(e)$. By the form of $h$ above, this value is odd if and only if $d$ and $e$ are of opposite parity. This completes the proof.
\end{proof}

\medskip

\begin{proof}[Proof of Theorem \ref{thm:classification_theorem}]
    The conjugacy classes of $\mathfrak{sl}_2$-triples in simple complex and real Lie algebras $\mathfrak{g}, \mathfrak{g}^\textbf{R}$ are studied via their one-to-one correspondence with nilpotent $G, G^\textbf{R}$-orbits \cite{Collingwood_McGovern_1993}. For classical complex Lie algebras $\mathfrak{g}=\mathfrak{sl}_n\mathrm{\textbf{C}},\,\mathfrak{so}_{2n+1\mathrm{\textbf{C}}},\, \mathfrak{sp}_{2n}\mathrm{\textbf{C}}, \,\mathfrak{so}_{2n}\mathrm{\textbf{C}}$ of type $A_n, B_n, C_n, D_n$ and their real forms $\mathfrak{g}^\textbf{R}$, nilpotent orbits are classified by Young diagrams associated to certain partitions of some $N$ and signed Young diagrams of signature $(p,q)$. A summary of the classification can be found in \cite[Theorem 3.3]{Bradlow_Collier_García-Prada_Gothen_Oliveira_2024}. We follow the convention used in \cite{Bradlow_Collier_García-Prada_Gothen_Oliveira_2024} that the Young diagram associated to the partition $N=\sum_i r_i\cdot i$ has $r_i$ rows of length $i$. For signed Young diagrams, let $p_i$ denote the number of rows of length $i$, whose leftmost box is labeled $+$, and $q_i$ the number of such rows whose leftmost box is labeled $-$. For exceptional Lie algebras $\mathfrak{g}=G_2,F_4, E_6,E_7,E_8$ and their real forms $\mathfrak{g}^\textbf{R}$, the nilpotent orbits are classified explicitly by their weighted Dynkin diagrams \cite{DOKOVIC1}, \cite{DOKOVIC2}. \\

    This proof essentially follows the classification of the magical triples, using the modified criteria in Proposition \ref{prop:modified_criteria}. The only significant difference is in the proof for exceptional Lie algebras, where the criteria is not as straightforwardly verified from the tables as in the magical case. Applying certain necessary conditions, we are able to rule out most of the triples from Tables VI-X in \cite{DOKOVIC1} and Tables VII-VIII in \cite{DOKOVIC2}, and for those that are left, we compute explicitly whether \eqref{eqn:extended_triple_criteria} holds. 
    
    We start with the classical case, namely, when $\hat{\rho}$ is an $\mathfrak{sl}_2$-triple of $\mathfrak{g}^\textbf{R}$ that is a real form of a classical Lie algebra $\mathfrak{g}$. We first impose the condition that $\mathfrak{c}^\textbf{R}$ is compact (the form of $\mathfrak{c}^\textbf{R}$ can be found in \cite[Section 3.2]{Bradlow_Collier_García-Prada_Gothen_Oliveira_2024}) and then check criteria \eqref{eqn:extended_triple_criteria}. 
    
    For $\mathfrak{g}^\textbf{R}=\mathfrak{sl}_n\mathrm{\textbf{R}}$, $\mathfrak{c}^\textbf{R}=\mathfrak{s}(\bigoplus_{i=1}^n \mathfrak{gl}_{r_i}\mathrm{\textbf{R}})$. So $\mathfrak{c}^\textbf{R}$ is compact if and only if the partition associated to $\hat{\rho}$ is $1\cdot n$, i.e. $r_n=1$ and $r_i=0$ for all $i\neq n$. This gives us the principal triples, which are even magical triples by Lemma \ref{lem:dimension_of_g0}.

    For $\mathfrak{g}^\textbf{R}=\mathfrak{su}^*_{2m}$, $\mathfrak{c}^\textbf{R}=\mathfrak{s}(\bigoplus_{i=1}^m \mathfrak{u}^*_{2r_i})$ is compact if and only if $r_m=1$ and $r_i=0$ for all $i\neq m$. By Lemma \ref{lem:dimension_of_g0}, this must be an even triple. Therefore, the classification follows from the classification of magical triples. 

    For $\mathfrak{g}^\textbf{R}=\mathfrak{su}(p,q)$, $q\ge p$, the centralizer $\mathfrak{c}^\textbf{R}=\mathfrak{s}(\bigoplus_{i=1}^{p+q} \mathfrak{u}(p_i,q_i))$ is compact if and only if for each $i$, $p_i=0$ or $q_i=0$. The $\mathfrak{sl}_2$-triple $\rho$ in $\mathfrak{sl}_{p+q}\mathrm{\textbf{C}}$ is associated to the partition $p+q=\sum_{i=1}^{p+q}r_i\cdot i$, where $r_i=p_i+q_i$. The dimension of the centralizer is 
    $\dim \mathfrak{c}=\sum_{i=1}^{p+q} r_i^2-1$ and $\dim \mathfrak{m}-\dim \mathfrak{h}=1-(q-p)^2$.
    Using the dimension of $\mathfrak{g}_0$ calculated in Lemma \ref{lem:dimension_of_g0}, the criteria Proposition \ref{prop:modified_criteria} asks that the following equality holds: 
    \begin{align*}
        1-(q-p)^2= \dim V_\rho-\sum_{\substack{i< j \\
        i\neq j \,(\mathrm{mod} \,2)}}\, 2i\,r_ir_j- 2\left(\sum_i r_i^2-1\right) 
    \end{align*}
    Rearranging the expression using the explicit formula of $\dim V_\rho$ from Theorem 6.1.3 of \cite{Collingwood_McGovern_1993}, we can rewrite the above equation as
    \begin{align*}
        -(q-p)^2=-r_1^2+\sum_{i\ge 3} (i-2)r_i^2+\sum_{\substack{i< j \\
        i= j \,(\mathrm{mod} \,2)}}\, 2i\,r_ir_j.
    \end{align*}
From this we conclude that either $r_1\neq 0$ or $q=p$ and $r_2$ is the only nonzero $r_i$. The latter corresponds to the even magical triple of $\mathfrak{su}(p,p)$. If $r_1\neq 0$, the nilpotent $\hat{e}$ must be an extended magical triple contained in an subalgebra isomorphic to $\mathfrak{su}(p,q-r_1)$. In this subalgebra, $\hat{e}$ has no $r_1$-term, so it must be that $q-r_1=p$ and $r_2=p$. Thus, $\hat{\rho}$ is an odd magical triple if and only if $\rho$ is associated to the partition $[2^p,1^{q-p}]$.

    For $\mathfrak{g}^\textbf{R}=\mathfrak{so}_{p,q}$, $\mathfrak{c}^\textbf{R}=\bigoplus_{i-even} \mathfrak{sp}_{p_i+q_i}\mathrm{\textbf{R}} \bigoplus_{i-odd} \mathfrak{so}_{p_i,q_i}$ is compact if and only if for all even $i$, $p_i=q_i=0$ and for each odd $i$, $p_i=0$ or $q_i=0$. By Lemma \ref{lem:dimension_of_g0}, this is always an even triple and the classification follows from the classification of magical triples.

    For $\mathfrak{g}^\textbf{R}=\mathfrak{so}^*_{2m}$, $\mathfrak{c}^\textbf{R}=\bigoplus_{i-even} \mathfrak{sp}_{2p_i,2q_i}\bigoplus_{i-odd} \mathfrak{so}^*_{2(p_i+q_i)}$ is compact if and only if for each even $i$, $p_i=0$ or $q_i=0$ and for odd $i$, $p_i+q_i \le 1$. The $\mathfrak{sl}_2$-triple $\rho$ in $\mathfrak{so}_{2m}\mathrm{\textbf{C}}$ is associated to the partition $2m=\sum_{i=1}^m (2r_i)\cdot i$, where $r_i=p_1+q_i$. The dimension of the centralizer is
    $$\dim \mathfrak{c}=\sum_{i=1}^m 2r_i^2+\sum_{i-even} r_i-\sum_{i-odd}r_i,$$
    and $\dim \mathfrak{m}-\dim \mathfrak{h}=-m$.
    Using the dimension $\mathfrak{g}_0$ calculated in Lemma 3.4, we check Criterion \eqref{eqn:extended_triple_criteria} in Proposition \ref{prop:modified_criteria}:
    \begin{align*}
-m=\dim V_\rho - \sum_{\substack{i< j \\
        i\neq j \,(\mathrm{mod} \,2)}}\, i\, r_ir_j-2\dim \mathfrak{c}.
    \end{align*}
    Using the explicit formula of $\dim V_\rho$ from Theorem 6.1.3 of \cite{Collingwood_McGovern_1993}, the above works out to be
    \begin{align*}
        -m=\sum_{i=1}^{p+q} \frac{1}{2}s_i^2-2\sum_{i=1}^{p+q} 2r_i^2+\sum_{i-odd} r_i-\sum_{i-even} 2r_i- \sum_{\substack{i< j \\
        i\neq j \,(\mathrm{mod} \,2)}}\, i\, r_ir_j.
    \end{align*}
    Rearranging the terms, we may rewrite the equation as 
    \begin{align*}
        -m=-2r_1^2+\sum_{i=3}^{p+q} 2(i-2)r_i^2+\sum_{i-odd} r_i- \sum_{i-even} 2r_i+\sum_{\substack{i< j \\
        i= j \,(\mathrm{mod} \,2)}}\, i\, r_ir_j.
    \end{align*}
    Since the expression $\sum_{i=3}^{p+q} 2(i-2)r_i^2-\sum_{i\ge 4,\,even} 2r_i \ge 0$, the right hand side is negative only if $r_1\neq 0$ or $r_2\neq 0$. If $r_1=0$, the right hand side of the expression is strictly positive if $r_i\neq 0$ for any $i\ge 3$. So the only nonzero $r_i$ is $r_2$, and we have $r_2=m/2$, which corresponds to the magical triple of $\mathfrak{so}^*_{2(m-r_1)}$, when $m$ is even.  If $r_1\neq 0$, the nilpotent $\hat{e}$ must be an extended magical triple in a subalgebra isomorphic to $\mathfrak{so}^*_{2m-2r_1}$, in which $\hat{e}$ has no $r_1$-term. With $r_1=0$, we are in the setting of the previous case, the only nonzero $r_i$ is $r_2$, and $r_2=(m-r_1)/2$. We solve for $r_1$ by the equality 
    \begin{align*}
        -m&=-2r_1^2+r_1-2r_2 \\
        &=-2r_1+r_1-(m-r_1),
    \end{align*}
    which works out to be $-r_1^2+r_1=0$. Hence, $r_1=0$ or $r_1=1$. Only the latter triple is odd. Thus, $\hat{\rho}$ is an odd magical triple if and only if $\rho$ is associated to the partition $[2^m,1^2]$.
    
    For $\mathfrak{g}^\textbf{R}=\mathfrak{sp}_{2m}\mathrm{\textbf{R}}$, $\mathfrak{c}^\textbf{R}=\bigoplus_{i-odd}\mathfrak{sp}_{p_i+q_i}\mathrm{\textbf{R}}\,
    \bigoplus_{i-even} \mathfrak{so}_{p_i,q_i}$ is compact if and only if for all odd $i$, $p_i=q_i=0$ and for each even $i$, $p_i=0$ or $q_i=0$. By Lemma \ref{lem:dimension_of_g0}, this is always an even triple  and the classification follows from the classification of magical triples.

    For $\mathfrak{g}^\textbf{R}=\mathfrak{sp}_{2p,2q}\mathrm{\textbf{R}}$, $\mathfrak{c}^\textbf{R}=\bigoplus_{i-odd}\mathfrak{sp}_{2p_i,2q_i}\mathrm{\textbf{R}} \bigoplus_{i-even} \mathfrak{so}^*_{2(p_i+q_i)}$ is compact if and only if for each odd $i$, $p_i=0$ or $q_i=0$ and for each even $i$, $p_i+q_i\le 1$. The $\mathfrak{sl}_2$-triple $\rho$ in $\mathfrak{sp}_{2m}\mathrm{\textbf{C}}$ is associated to the partition $2p+2q=\sum_{i=1}^{p+q} 2r_i\cdot i$, where each $r_i=p_i+q_i$. The dimension of the centralizer is
    $$\dim \mathfrak{c}=\sum_{i=1}^{p+q}2r_i^2-\sum_{i-even}r_i+\sum_{i-odd} r_i. $$
    Using the dimension $\mathfrak{g}_0$ calculated in Lemma \ref{lem:dimension_of_g0}, we check the criteria in Proposition \ref{prop:modified_criteria}:
    \begin{align*}
        -p-q-2(q-p)^2=\dim V_\rho - \sum_{\substack{i< j \\
        i\neq j \,(\mathrm{mod} \,2)}}\, i\, r_ir_j-2\dim \mathfrak{c}.
    \end{align*}
    Using the explicit formula of $\dim V_\rho$ from Theorem 6.1.3 of \cite{Collingwood_McGovern_1993}, the above is equivalent to 
    \begin{align*}
         -p-q-2(q-p)^2=\sum_{i=1}^{p+q} \frac{1}{2}s_i^2-2\sum_{i=1}^{p+q} 2r_i^2-\sum_{i-odd} r_i+\sum_{i-even} 2r_i- \sum_{\substack{i< j \\
        i\neq j \,(\mathrm{mod} \,2)}}\, i\, r_ir_j.
    \end{align*}
    Rearranging the terms, we may rewrite the equation as 
    \begin{align*}
         -p-q-2(q-p)^2=-2r_1^2+\sum_{i=3}^{p+q} 2(i-2)r_i^2-\sum_{i-odd} r_i+ \sum_{i-even} 2r_i+\sum_{\substack{i< j \\
        i= j \,(\mathrm{mod} \,2)}}\, i\, r_ir_j.
    \end{align*}
    Since the expression $\sum_{i=3}^{p+q} 2(i-2)r_i^2-\sum_{i\ge 3,\,odd} r_i \ge 0$, the right hand side is negative only if $r_1\neq 0$. Thus the nilpotent $\hat{e}$ is be an extended magical triple in a subalgebra isomorphic to $\mathfrak{sp}_{2p-2r_1,2q}$ or $\mathfrak{sp}_{2p,2q-2r_1}$, with $r_1=0$. But the previous argument says that there are no extended magical triples for $r_1\neq 0$. We therefore conclude that there are no extended magical triples in this case. This concludes the classification of odd magical triples in classical Lie algebras. The weighted Dynkin diagram of the triples can be computed using the algorithm from \cite[Section 5.3]{Collingwood_McGovern_1993}. \\
    
    For the exceptional case, we consult Tables VI-XV of \cite{DOKOVIC1} for real forms of inner type and Tables VII and VIII of \cite{DOKOVIC2} for the two real forms of outer type: $E_6^6$ and $E_6^{-26}$. In the former case, given an $\mathfrak{sl}_2$-triple $\hat{\rho}$ of a real form $\mathfrak{g}^\textbf{R}$ of exceptional Lie algebra $\mathfrak{g}=G_2,F_4,E_6,E_7,E_8$, the superscripts on the Lie algebra is the difference $\dim \mathfrak{m}-\dim \mathfrak{h}$, Column 1 lists the weighted Dynkin diagram associated to $\rho$, Column 4 column lists $\dim V_\rho \cap \mathfrak{h}$, Column 5 lists $\dim\mathfrak{c}\cap \mathfrak{h}$ and Column 6 lists $\dim V_{even}\cap \mathfrak{m}$ and the last column lists the isomorphism types of $\mathfrak{c}^\textbf{R}\cap \mathfrak{h}^\textbf{R}$ (double check this last one). An $\mathfrak{sl}_2$-triple is extended magical if and only if it satisfies the following conditions:
    \begin{itemize}
        \item [(a)] Column 5 is equal to $\dim \mathfrak{c}$ and the last column is a compact Lie algebra,
        \item [(b)] Column 6 is equal to $\dim V_{even}$, i.e. $\dim V_{even}=\dim V_{even}\cap \mathfrak{m}$,
        \item [(c)] The difference $(\text{Column 6} - \text{Column 5})$ is equal to the superscript of the Lie algebra.
    \end{itemize}
     Notice that, assuming (a) and (b) are satisfied, condition (c) is equivalent to \eqref{eqn:extended_triple_criteria}. Since the even magical triples have already been classified, we search among the odd triples in Tables VI-XV, first imposing conditions (a) and (c). The only possible odd magical triples are listed below, along with their $\mathfrak{sl}_2$-data, which are computed using their root poset diagrams:
     \begin{center}
\begin{tabular}{ |c|c|c|c| } 
  \hline $\mathfrak{g}$ & $E_6$ & 
  $E_7$ & $E_8$\\ 
  \hline $\mathfrak{g}^\textbf{R}$ & $E_6^{-14}$ & 
  $E_7^7$ & $E_8^8$\\
  \hline
  Table & X & XI & XIV \\
  \hline 
  Rows & 3, 4 & 48, 49 & 42 \\
  \hline 
  $\rho$ & 10001 & 1001010 & 00010001 \\ 
  \hline $\mathfrak{c}^\textbf{R}$ & $\mathfrak{so}_7\mathrm{\textbf{R}}+\mathfrak{so}_2\mathrm{\textbf{R}}$ & $\mathfrak{so}_2\mathrm{\textbf{R}}+\mathfrak{so}_2\mathrm{\textbf{R}}$ & $\mathfrak{su}(2)+\mathfrak{so}_2\mathrm{\textbf{R}}$ \Tstrut\Bstrut \\
  \hline\hline
  $n_0$ & 22 & 2 & 4 \\
  \hline
  $n_1$ & 16 & 5 & 8\\
  \hline
  $n_2$ & 8 & 3 & 9\\
  \hline
  $n_3$ & 0 & 3 & 8\\
  \hline
  $n_4$ & 0 & 6 & 9\\
  \hline
  $n_5$ & 0 & 4 & 8\\
  \hline
  $n_6$ & 0 & 3 & 5\\
  \hline
  $n_7$ & 0 & 2 & 4\\
  \hline
  $n_8$ & 0 & 1 & 1\\
  \hline 
\end{tabular}
\end{center}
Imposing condition (b) and comparing Column 6 to the dimension of the nontrivial even highest weight spaces, only the triple of $e^{-14}_6$ turns out to be odd magical. \\

For the two real forms of outer type: $E_6^6$ and $E_6^{26}$, we consult Tables VII and VII of \cite{DOKOVIC2} for relevant information. The first column lists the weighted Dynkin diagram associated to the triple $\rho$, the third to last column lists $\dim V_\rho \cap \mathfrak{h}$ and the last column lists the isomorphism type of $\mathfrak{c}^\textbf{R}$. Observe that apart from the known even magical triple listed in Row 20 of Table VIII, $E_6^6$ does not admit any odd triples whose $\mathfrak{c}^\textbf{R}$ is compact. On the other hand $E_6^{-26}$ has no even magical triples and only one odd triple, listed in Row 1 of Table VII. This triple shares the same complex weighted Dynkin diagram as the odd magical triple in $E^{-14}_6$. However, from the $\mathfrak{sl}_2$-data above, we see that $n_0=22\neq \dim \mathfrak{c}^\textbf{R}\cap \mathfrak{h}^\textbf{R}=21$. Hence condition (a) is not satisfied and this triple is not an odd magical triple. This concludes the classification of odd magical triples in exceptional Lie algebras. 
\end{proof}

\begin{remark}
    Each of the odd magical triples is induced by trivially extending the magical $\mathfrak{sl}_2$-triple of its corresponding Hermitian Lie algebra of tube type. This can be observed from the proof of the classification and can be clearly seen from the explicit $\mathfrak{sl}_2$-data in the next section. Conversely, one can also deduce from the classification that the three odd magical triples are the only $\mathfrak{sl}_2$-triples retaining interesting properties by extending the even magical triples trivially.
\end{remark} 

\begin{remark}
Folding with respect to the involution defining $\mathfrak{g}^\textbf{R}$, the weighted Dynkin diagrams of the restricted root system of the odd magical triples are given by the non-reduced root systems:
\begin{align*}
\mathfrak{g}^\textbf{R}=\mathfrak{su}(p,q), \quad &BC_{p}: \dynkin [labels*={,,,,\alpha}]C{} \\
\mathfrak{g}^\textbf{R}=\mathfrak{so}^*_{4n+2}, \quad &BC_n: \dynkin [labels*={,,,,\alpha}]C{} \\
\mathfrak{g}^\textbf{R}=E_6^{-14}, \quad &BC_2: \dynkin[labels*={,\alpha}] C2
\end{align*}
Here the simple restricted root corresponding to the simple complex root of weight-1 is the longest root $\alpha$. From Remark \ref{rmk:relation_to_theta_positivity}, we can identify $2\alpha$ as the restricted root whose associated unipotent subalgebra $\mathfrak{u}_{2\alpha}$ contains an acute convex cone invariant under the action of $L^0_\alpha$ \cite[Section 3.7]{Guichard_Wienhard_2025}. This phenomenon will be investigated in forthcoming work. 
\end{remark}

\section{Odd magical triples}
\label{sec:odd_magical_triples}
We calculate the explicit $\mathfrak{sl}_2$-data of the odd magical triples and compare them to the $\mathfrak{sl}_2$-data of the even magical triple of its maximal subtube.

\begin{proposition}
    The $\mathfrak{sl}_2$-data of an odd magical triple of $\mathfrak{g}^\textbf{R}$ along with the $\mathfrak{sl}_2$-data of their maximal subtube $\mathfrak{g}_t^\textbf{R}$ are listed below. The last row lists the number $s=\dim \mathfrak{m}-\dim \mathfrak{h}$, which could be found in Table 2 of \cite{Bradlow_Collier_García-Prada_Gothen_Oliveira_2024}. Noting that $n_0=\dim \mathfrak{c}$ and $n_2=\dim \mathfrak{g}_0-\dim \mathfrak{c}$, Criterion \eqref{eqn:extended_triple_criteria} can be easily verified for all the cases below by observing that in each case $s=n_2-n_0$.
    \begin{enumerate}
        \item[(1)] $\mathfrak{g}^\textbf{R}=\mathfrak{su}(p,q)$, $q> p$, $\mathfrak{g}=\mathfrak{sl}_{p+q}\textbf{C}$:

\begin{center}
\begin{tabular}{ |c|c|c| } 
\hline & $\mathfrak{g}_t^\textbf{R}$ & $\mathfrak{g}^\textbf{R}$ \\
  \hline 
  & $\mathfrak{su}(p,p)$ & 
  $\mathfrak{su}(p,q)$ \\ 
  \hline $\otimes\mathrm{\textbf{C}}$ & $\mathfrak{sl}_{2p}\mathrm{\textbf{C}}$ & 
  $\mathfrak{sl}_{p+q}\mathrm{\textbf{C}}$ \\ 
  \hline \hline
  $\rho$ & $[2^p]$ & $[2^p,1^{q-p}]$ \\
  \hline
  $\mathfrak{c}$ & $\mathfrak{sl}_p\mathrm{\textbf{C}}$ & $\mathfrak{s}(\mathfrak{gl}_p\mathrm{\textbf{C}}\oplus \mathfrak{gl}_{q-p}\mathrm{\textbf{C}})$  \\
  \hline
$\mathfrak{c}^\textbf{R}$ & $\mathfrak{su}(p)$ & $\mathfrak{s}(\mathfrak{u}(p)\oplus \mathfrak{u}(q-p))$  \\
  \hline \hline
  $n_0$ & $p^2-1$ & $p^2-1+(q-p)^2$ \\
  \hline
  $n_1$ & $0$ & $2p(q-p)$ \\
  \hline
  $n_2$ & $p^2$ & $p^2$ \\
  \hline
  $s$ & $1$ & $1-(q-p)^2$  \\
  \hline
\end{tabular}

\end{center}
\item [(2)] $\mathfrak{g}^\textbf{R}=\mathfrak{so}^*_{4n+2}$, $\mathfrak{g}=\mathfrak{so}^*_{4n}$:
\begin{center}
\begin{tabular}{ |c|c|c| } 
\hline & $\mathfrak{g}_t^\textbf{R}$ & $\mathfrak{g}^\textbf{R}$ \\
\hline  & $\mathfrak{so}^*_{4n}$ & $\mathfrak{so}^*_{4n+2}$ \\ 
  \hline $\otimes \mathrm{\textbf{C}}$ & $\mathfrak{so}_{4n}\mathrm{\textbf{C}}$ & 
  $\mathfrak{so}_{4n+2}\mathrm{\textbf{C}}$ \\ 
  \hline \hline
  $\rho$ & $[2^{2n}]$ & $[2^{2n},1^{2}]$ \\ 
  \hline
  $\mathfrak{c}$ & $\mathfrak{sp}_{2n}\mathrm{\textbf{C}}$ & $\mathfrak{sp}_{2n}\mathrm{\textbf{C}} \oplus \mathrm{\textbf{C}}$ \\ 
  \hline
    $\mathfrak{c}^\textbf{R}$ & $\mathfrak{sp}_{2n}\mathrm{\textbf{R}}$ & $\mathfrak{sp}_{2n}\mathrm{\textbf{R}} \oplus \mathfrak{u}(1)$ \\ 
  \hline \hline
  $n_0$ & $n(2n+1)$ & $n(2n+1)+1$ \\
  \hline
  $n_1$ & $0$ & $4n$ \\
  \hline
  $n_2$ & $n(2n-1)$ & $n(2n-1)$ \\
  \hline 
  $s$ & $-2n$ & $-2n-1$ \\
  \hline
\end{tabular}
\end{center}

\newpage
\item [(3)]
$\mathfrak{g}^\textbf{R}=E_6^{-14}$, $\mathfrak{g}=E_6$:
\begin{center}
\begin{tabular}{ |c|c|c | } 
\hline & $\mathfrak{g}_t^\textbf{R}$ & $\mathfrak{g}^\textbf{R}$ \\ 
  \hline  & $\mathfrak{so}_{2,8}$ & 
  $E_6^{-14}$ \\ 
  \hline $\otimes \mathrm{\textbf{C}}$ & $\mathfrak{so}_{10}\mathrm{\textbf{C}}$ & 
  $E_6$ \\
  \hline \hline
  $\rho$ & $[3^1,1^7]$ & $100001$ \\ 
  \hline
  $\mathfrak{c}$ & $\mathfrak{so}_7\mathrm{\textbf{C}}$ & $\mathfrak{so}_{7}\mathrm{\textbf{C}} \oplus \mathrm{\textbf{C}}$ \\ 
  \hline
$\mathfrak{c}^\textbf{R}$ & $\mathfrak{so}_7\mathrm{\textbf{R}}$ & $\mathfrak{so}_7\mathrm{\textbf{R}} \oplus \mathfrak{u}(1)$ \\
  \hline \hline
  $n_0$ & $21$ & $22$ \\
  \hline
  $n_1$ & $0$ & $16$ \\
  \hline
  $n_2$ & $8$ & $8$ \\
  \hline 
  $s$ & $-13$ & $-14$ \\
  \hline
\end{tabular}
\end{center}
    \end{enumerate}
\end{proposition}
\smallskip

Recall that the normalizer $\mathfrak{g}_T^\textbf{R}$ of the maximal subtube $\mathfrak{g}_t^\textbf{R}$ is the maximal Hermitian Lie subalgebra of $\mathfrak{g}^\textbf{R}$. Over the $\mathrm{\textbf{Z}}$-grading of $\mathfrak{g}$, notice that the complexification of $\mathfrak{g}_T^\textbf{R}$ is simply the even $\mathrm{ad}_h$-weight spaces $\mathfrak{g}_T=N_G(\mathfrak{g_t})=\mathfrak{g}_{-2}\oplus\mathfrak{g}_0\oplus \mathfrak{g}_2$, and an odd magical triple of $\mathfrak{g}^\textbf{R}$ is an even magical triple of $\mathfrak{g}_T^\textbf{R}$ upon restriction. Therefore, we have the following proposition on the Lie group level:

\begin{proposition}
\label{prop:restriction_of_odd_magical}
    An odd magical triple of $G^\textbf{R}$ restricts to an even magical triple of the normalizer $G_T^\textbf{R}$ of its maximal subtube $G_t^\textbf{R}$.
\end{proposition}

The property of being extended magical give rise to the Lie algebra involution $\theta:\mathfrak{g}_0\to \mathfrak{g}_0$ given by \eqref{eq:theta_involution}.
Let $\tilde{\mathfrak{g}}^\textbf{R}$ be the semisimple part of the Cayley real form $\mathfrak{g}_{\mathcal{C}}^\textbf{R}$ and on the level of Lie groups, let $\tilde{G}^\textbf{R}$ be the subgroup of $G^\textbf{R}$ with Lie algebra $\tilde{\mathfrak{g}}^\textbf{R}$ and maximal compact subgroup $\mathcal{C}:=C\cap G^\textbf{R}$. Since $\mathfrak{g}^\textbf{R}_\mathcal{C}\subset \mathfrak{g}_0\subset \mathfrak{g}_T\subset \mathfrak{g}$, we have that $\tilde{\mathfrak{g}}_T^\textbf{R}=\tilde{\mathfrak{g}}^\textbf{R}$ and $\tilde{G}^\textbf{R}_T=\tilde{G}^\textbf{R}$.

From the $\mathfrak{sl}_2$-data above, observe that $\mathfrak{c}=\mathfrak{c}_t\oplus \hat{\mathfrak{c}}$, where $\hat{\mathfrak{c}}=\mathfrak{h}_T\cap \ker(\mathrm{ad})|_{\mathfrak{m}_t}$ those elements acting trivially on $\mathfrak{m}_t$, from which we have that $\mathfrak{h}_T=\mathfrak{h}\oplus \hat{\mathfrak{c}}$, and $\mathfrak{g}_T=\mathfrak{g}_t\oplus \hat{\mathfrak{c}}.$ On the level of Lie groups, let $\hat{C}:=H_T\cap \ker(\mathrm{Ad})|_{\mathfrak{m}_t}$, then we have that 
\begin{align*}
    H_T=H_t\times_Z  \hat{C}, \\
    G_T=G_t\times_{Z}  \hat{C},
\end{align*}
where the quotient is taken with respect to the diagonal action of the center $Z:=Z(G)$ restricted to their common central subgroup. Let $\hat{\mathfrak{c}}^\textbf{R}$ be the real form of $\hat{\mathfrak{c}}$ associated to $\sigma$, then we have that 
\begin{align*}
\mathfrak{h}_T^\textbf{R}=\mathfrak{h}^\textbf{R}_t\oplus \hat{\mathfrak{c}}^\textbf{R}, \\
\mathfrak{g}_T^\textbf{R}=\mathfrak{g}^\textbf{R}_t\oplus  \hat{\mathfrak{c}}^\textbf{R},\\
    \tilde{\mathfrak{g}}^\textbf{R}=\tilde{\mathfrak{g}}_t^\textbf{R}\oplus \hat{\mathfrak{c}}^\textbf{R}.
\end{align*}
We list the data of $\tilde{\mathfrak{g}}^\textbf{R}$ and $\hat{\mathfrak{c}}^\textbf{R}$ of each odd magical triple below:
\begin{table}[htb]
\centering
\begin{tabular}{ |c|c|c|c| } 
\hline $\mathfrak{g}^\textbf{R}$ & $\mathfrak{g}^\textbf{R}_t$ & $\tilde{\mathfrak{g}}^\textbf{R}$ & $\hat{\mathfrak{c}}^\textbf{R}$ \Tstrut\Bstrut\\
  \hline \hline $\mathfrak{su}(p,q)$ & $\mathfrak{su}(p,p)$ & $\mathfrak{sl}_p\mathrm{\textbf{C}}\oplus \hat{\mathfrak{c}}^\textbf{R} $ & $\mathfrak{s}(\mathfrak{u}(p-q)\oplus\mathfrak{u}(1)) $ \Tstrut \Bstrut\\
\hline $\mathfrak{so}^*_{4m+2}$ & $\mathfrak{so}^*_{4m}$ & $\mathfrak{su}^*_{2m}\oplus \hat{\mathfrak{c}}^\textbf{R}$ & $\mathfrak{u}(1)$ \Tstrut \Bstrut \\
\hline 
$\mathrm{E}_6^{-14}$ & $\mathfrak{so}_{2,8}$ & $\mathfrak{so}_{1,7}\oplus\hat{\mathfrak{c}}^\textbf{R}$ & $\mathfrak{u}(1)$ \Tstrut\Bstrut \\
\hline
\end{tabular}
\caption{Subalgebras associated to odd magical triples of $\mathfrak{g}^\textbf{R}$}
\label{table:data_of_odd_magical_triples}
\end{table}

\section{Geometric characterization of extended magical triples}
In this section we give a characterization of the extended magical triples via intersections of their associated Slodowy slice with the $G^\textbf{R}$-Higgs bundles moduli spaces. Namely, we show that, given that the underlying curve has sufficiently large genus of the underlying curve, an $\mathfrak{sl}_2$-triple is an extended magical triple if and only if its associated Slodowy slice $\mathrm{Slo}_\rho$ is a union of connected components in $\mathcal{M}(G^\textbf{R})$. We begin by proving the following statement: 

\begin{theorem}\label{thm:maximal_components}
    Let $G$ be a simple complex Lie group and let $\hat{\rho}$ be an odd magical triple of the real form $G^\textbf{R}\subset G$ with $\rho$ its Cayley transform. The Slodowy slice $\mathrm{Slo}_\rho \subset \mathcal{M}(G^\textbf{R})$ is precisely the union of maximal $G^\textbf{R}$-Higgs bundles $\mathcal{M}_{\mathrm{max}}(G^\textbf{R})$.
\end{theorem}

\begin{lemma}
\label{lem:intersection}
   Let $\hat{\rho}$ be an odd magical triple of $G^\textbf{R}$. The image of the Slodowy map
\begin{align*}
\hat{\Psi}_\rho:\,\,&\mathcal{B}_\rho(G) \to \mathcal{H}(G), \\ &((E_C,\varphi_C),  \{\varphi_{j,\,n_j}\}_{j=1}) \mapsto (E_G:=(E_C \times E_T)[G], \underline{f}+\varphi_C+\sum_{(j,n_j)} \varphi_{j,\,n_j}),
\end{align*}
is in $\mathcal{H}(G^\textbf{R})$ if and only if $E_C$ reduces to a $C\cap H$-bundle, $\varphi_C=0$, and each holomorphic section $\varphi_{j,n_j} \in H^0(X, E_C[{v_j}^{n_j}]\otimes K)$ that is not valued in $\mathfrak{m}$ vanishes. 
\end{lemma}

\begin{proof}
    This is a straightforward observation using the definition of a $G^\textbf{R}$-Higgs pair that $E_G=(E_C\times E_T)[G]$ must reduce to a holomorphic $H$-bundle $E_H$ and the Higgs field $f+\varphi_C+\sum_{(j,n_j)} \varphi_{j,n_j}$ must be a holomorphic section in $H^0(X, E_H[\mathfrak{m}]\otimes K)$. Since $T \subset H$, the first condition is equivalent to asking that
    $$E_C\cong E_{C\cap H}[C],$$
    i.e. $E_C$ reduces to a $C\cap H$-bundle. The condition on the Higgs field requires that $\varphi_C\in H^0(X,E_C[\mathfrak{c}]\otimes K)$ vanish, as $\mathfrak{c}\cap \mathfrak{m}=0$. The same applies to all the other sections $\varphi_{j,n_j}$ not valued in $\mathfrak{m}$.
\end{proof}

The proof of Theorem \ref{thm:maximal_components} can be broken down into the following two propositions.

\begin{proposition}
\label{prop:part1}
Let $\hat{\rho}$ be an odd magical triple of $G^\textbf{R}$. Consider the Slodowy map restricted to the preimage of polystable $G^\textbf{R}$-Higgs pairs in $\mathcal{H}(G^\textbf{R})$ and descending to the moduli space $\mathcal{M}(G^\textbf{R})$. Every point in the image $\mathrm{Slo}_\rho\subset \mathcal{M}(G^\textbf{R})$ of the restricted Slodowy map has maximal Toledo number. \end{proposition}

\begin{proof}
By the previous Lemma, the image of the restricted Slodowy map in $\mathcal{H}(G^\textbf{R})$ consist of elements of the form
\begin{align*}
((E_C \times E_T)[H], \underline{f}+\sum_{j=1}^N \varphi_{j,n_j}).
\end{align*}
Let $\mathcal{E}_{G^\textbf{R}}=((E_C \times E_T)[H], \underline{f}+\sum_{j=1}^N \varphi_{j,n_j})$ be a polystable $G^\textbf{R}$-Higgs bundle in the Slodowy slice $\mathrm{Slo}_\rho\subset \mathcal{M}(G^\textbf{R})$. Since $C=C_t\times_{Z} \hat{C}$, there exist a $C_t$-bundle and a $\hat{C}$-bundle such that $E_C=E_{C_t}\times_{Z}E_{\hat{C}}$. Now $C$ commutes with $T$, thus
\begin{align*}
(E_{C_t}\times_{Z}E_{\hat{C}})\times E_T=(E_{C_t}\times E_T)\times_{Z} E_{\hat{C}}.
\end{align*}

For $\mathcal{E}_{G_t^\textbf{R}}=((E_{C_t} \times E_T)[H], \underline{f}+\sum_{j=1}^N \varphi_{j,n_j})$ and $\mathcal{E}_{\hat{C}^\textbf{R}}=(E_{\hat{C}},0)$, we may write 
$$\mathcal{E}_{G^\textbf{R}}=\mathcal{E}_{G_t^\textbf{R}}\times_{Z} \mathcal{E}_{\hat{C}^\textbf{R}}.$$

The Toledo character of $\mathcal{E}_{G^\textbf{R}}$ breaks into the sum of the Toledo characters of $\mathcal{E}_{G_t^\textbf{R}}$ and  $\mathcal{E}_{\hat{C}^\textbf{R}}$, $$\chi=\chi_t+\chi_{\hat{C}}.$$  
Notice that the central subgroup $Z=Z(G_T)$ is in the kernel of the lifted Toledo characters $\tilde{\chi}_t$ and $\tilde{\chi}_{\hat{C}}$, and the lifted Toledo character $\tilde{\chi}:H\to \textbf{C}$ is given by the product of their respective lifts $$\tilde{\chi}=\tilde{\chi}_t\cdot \tilde{\chi}_{\hat{C}}:H_t\times_{Z} \hat{C}\to \textbf{C}.$$ 

Therefore, the Toledo number can be computed as 
\begin{align*}
\tau(\mathcal{E}_{G^\textbf{R}})=\deg((E_C\times E_T)(\tilde{\chi}))=\deg((E_{C_t}\times E_T)(\tilde{\chi}_t))+\deg(E_{\hat{C}}(\tilde{\chi}_{\hat{C}}))=\tau(\mathcal{E}_{G_t^\textbf{R}})+\tau(\mathcal{E}_{\hat{C}^\textbf{R}}).
\end{align*}
The first summand achieves the maximal Toledo number since $\mathcal{E}_{G_t^\textbf{R}}$ is a maximal $G_t^\textbf{R}$-Higgs bundle by the Cayley correspondence of Hermitian Lie groups in Example \ref{ex:Cayley_correspondence_Hermitian}, and the maximal Toledo number for a $G^\textbf{R}$-Higgs bundle of a nontube type Hermitian Lie group is the same as the maximal Toledo number of that of its maximal subtube. Furthermore, by maximality of the first summand, the second summand must be zero, for otherwise the above expression would violate the Milnor--Wood inequality by considering the first summand at the other extremal value. This shows that $\tau(\mathcal{E}_{G^\textbf{R}})$ is a maximal $G^\textbf{R}$-Higgs bundles. Observe that the vanishing of $\tau(\mathcal{E}_{\hat{C}^\textbf{R}})$ agrees with the statement in Theorem \ref{thm:maximal_Higgs_bundles}.
\end{proof}

\begin{proposition}
\label{prop:part2}
Any maximal polystable $G^\textbf{R}$-Higgs bundle $\mathcal{E}_{G^\textbf{R}}$ is in the image $\mathrm{Slo}_\rho$ of the restricted Slodowy map.
\end{proposition}

\begin{proof}
Recall from Theorem \ref{thm:maximal_Higgs_bundles} that any polystable maximal $G^\textbf{R}$-Higgs bundle for a nontube type Hermitian Lie group reduces to a $G_T^\textbf{R}$-Higgs bundle. Using the decomposition $H_T=H_t\times_{Z} \hat{C}$, we have that that for any maximal $G^\textbf{R}$-Higgs bundle $\mathcal{E}_{G^\textbf{R}}=(E_H,\varphi_H)$, there exist a maximal $G_t^\textbf{R}$-Higgs bundle $\mathcal{E}_{G_t^\textbf{R}}=(E_{H_t},\varphi_{H_t})$ and a $\hat{C}^\textbf{R}$-Higgs bundle $\mathcal{E}_{\hat{C}^\textbf{R}}=(E_{\hat{C}},\varphi_{\hat{C}})$ such that 
    \begin{align}
    \label{eqn:sum_of_Higgs_bundles}
\mathcal{E}_{G^\textbf{R}}=\mathcal{E}_{G_t^\textbf{R}}\times_{Z}\mathcal{E}_{\hat{C}^\textbf{R}}=(E_{H_t}\times_{Z} E_{\hat{C}},\varphi_{H_t}+\varphi_{\hat{C}})
    \end{align}
where $E_{H_t}\times_{Z} E_{\hat{C}}$ is a holomorphic $G_T^\textbf{R}$-bundle and $\varphi_{H_t}+\varphi_{\hat{C}}\in H^0(X,(E_{H_t}\times_Z E_{\hat{C}})[\mathfrak{m}_t]\otimes K)$. 
Since $\hat{\mathfrak{c}}^\textbf{R}\subset \mathfrak{c}^\textbf{R}\subset \mathfrak{h}^\textbf{R}$ is a compact Lie algebra, $\hat{C}^\textbf{R}$ is compact for any odd magical triple, so the Higgs field of $\mathcal{E}_{\hat{C}^\textbf{R}}$ vanishes $\varphi_{\hat{C}}=0$. By the Cayley correspondence for maximal $G^\textbf{R}_t$-Higgs bundles, we have that
\begin{align*}
\mathcal{E}_{G^\textbf{R}_t}=((E_{C_t}\times E_T)[H_t], \underline{f}+\sum_{j=1}^N\varphi_{j,n_j}).
\end{align*}

Therefore, \eqref{eqn:sum_of_Higgs_bundles} can be rewritten as
\begin{align*}
\mathcal{E}_{G^\textbf{R}}=((E_{C_t}\times E_T)\times_Z E_{\hat{C}})[H],\underline{f}+\sum_{j=1}^N\varphi_{j,n_j}),
\end{align*}
which is an element in the Slodowy slice $\mathrm{Slo}_\rho$ by Lemma \ref{lem:intersection}.
\end{proof}

\begin{proof}[Proof of Theorem \ref{thm:maximal_components}] The Theorem is proven by combining Propositions \ref{prop:part1} and \ref{prop:part2}.
\end{proof}

\medskip

We now shift our perspective by fixing real form $G^\textbf{R}\subset G$ and considering all of its $\mathfrak{sl}_2$-triples. If an $\mathfrak{sl}_2$-triple of $G^\textbf{R}$ is extended magical, then by the Cayley correspondence and Theorem \ref{thm:maximal_components}, its associated Slodowy slice intersects $\mathcal{M}(G^\textbf{R})$ at a union of connected components. It turns out that, with an assumption on the genus of the underlying curve, the converse is also true, and we have the following characterization of extended magical triples:
\begin{theorem}
\label{thm:characterization}
    Let $G^\textbf{R} \subset G$ be a real form of a complex simple Lie group with Lie algebra $\mathfrak{g}^\textbf{R}$. If an $\mathfrak{sl}_2$-triple $\hat{\rho}$ of $\mathfrak{g}^\textbf{R}$ is extended magical, then its Slodowy slice $\mathrm{Slo}_\rho \subset \mathcal{M}(G^\textbf{R})$ is a union of connected components. Conversely, assuming that $X$ has genus $g(X)\ge 2\dim_\textbf{R}(G^\textbf{R})^2$, then $\mathrm{Slo}_\rho \subset \mathcal{M}(G^\textbf{R})$ is a union of connected components only if the $\mathfrak{sl}_2$-triple $\hat{\rho}$ of $\mathfrak{g}^\textbf{R}$ is extended magical.
\end{theorem}

\begin{remark}
Notice that the trivial triple, which is excluded from the discussion in \cite{Bradlow_Collier_García-Prada_Gothen_Oliveira_2024}, is an even magical triple for trivial reasons. In view of Theorem \ref{thm:characterization}, the Slodowy slice of the trivial triple $\mathrm{Slo}_{\rho_{\mathrm{triv}}}$ is the connected component of $H^\textbf{R}$-Higgs bundles $\mathcal{N}(H^\textbf{R})\subset \mathcal{M}(G^\textbf{R})$ consisting of Higgs pairs $(E_H,0)$ with vanishing Higgs fields, which is, however, not a higher Teichm\"uller component.
\end{remark} 

In the following proposition, we show that the Slodowy slice $\mathrm{Slo}_\rho$ agrees with the expected dimension of $\mathcal{M}(G^\textbf{R})$ precisely when $\hat{\rho}$ is an even magical triple. 
\begin{proposition}
\label{prop:even_magical_triples}
 Let $G^\textbf{R}\subset G$ be the real form of a simple complex Lie group with Lie algebra $\mathfrak{g}^\textbf{R}$. An $\mathfrak{sl}_2$-triple $\hat{\rho}$ is an even magical triple of $\mathfrak{g}^\textbf{R}$ if and only if the Slodowy slice $\mathrm{Slo}_\rho \subset \mathcal{M}(G^\textbf{R})$ is of real dimension $2(g-1)\dim_\textbf{R} G^\textbf{R}$.
\end{proposition}

\begin{proof}
Suppose
$\hat{\rho}$ is not an even magical triple of $\mathfrak{g}^\textbf{R}$, we show that the dimension of $\mathrm{Slo}_\rho \subset  \mathcal{M}(G^\textbf{R})$ is strictly less than $2(g-1)\dim \mathfrak{g}^\textbf{R}$. By Lemma \ref{lem:intersection}, the image of the Slodowy map in $\mathcal{M}(G^\textbf{R})$ consists of elements of the form
\begin{align*}
((E_C \times E_T)[H], \underline{f}+\sum_{j=1}^N \varphi_{j,n_j}),
\end{align*}
where $E_C$ reduces to a $C\cap H$-bundle and each $\varphi_{j,n_j}\in H^0(X,E_C[{v_j}^{n_j}]\otimes K)$ is valued in $\mathfrak{m}$. 

Let $A=\bigoplus_{j=1}^{N} A_j\subset \mathfrak{m}\cap V_\rho^{>0}$ be the largest subspace preserved by the action of $C$, where $A_j$ is spanned by the eigenbasis $A_j:=A\cap V_j=\textrm{span}({v_j}^{n_{j_i}})$. For each ${v_j}^{n_{j_i}}\in A_j$, $H^0(X, E_C[{v_j}^{n_{j_i}}]\otimes K)\cong H^0(K^{j+1})$, whose dimension can be calculated via Riemann--Roch as
$$\dim H^0(K^{j+1})=2(g-1)(2j+1).$$

The dimension of $\mathrm{Slo}_\rho$ can be calculated and compared to the expected dimension by separating the contributions from the bundle part and from the Higgs field part: 
\begin{align*}
    \dim (\mathrm{Slo}_\rho)&=2(g-1)\dim_\textbf{C}(C\cap H)+2(g-1)\sum_{j=1}^N \dim_\textbf{C} A_j \cdot (2j+1) \\
&< 2(g-1)\dim_\textbf{C}C+2(g-1)\sum_{j=1}^N\dim_\textbf{C}V_{j}\cdot (2j+1)\\
&=2(g-1)\dim_{\textbf{C}}G
\\
&=2(g-1)\dim_{\textbf{R}}G^\textbf{R},
\end{align*}
where the inequality is strict because by our assumption that $\hat{\rho}$ is not an even magical triple, it must be that either $f+V_{\rho} \not\subset \mathfrak{m}$ or $\mathfrak{c} \not\subset\mathfrak{h}$, equivalently, either $\dim_\textbf{C} A<\dim_\textbf{C} V_\rho^{>0}$ or $\dim_\textbf{C} C\cap H < \dim_\textbf{C} C$. 

In the other direction, suppose $\hat{\rho}$ is an even magical triple of $\mathfrak{g}^\textbf{R}$, then $A=V_\rho^{>0}$ and $\dim_\textbf{C} C\cap H=\dim_\textbf{C} C$ and the above strict inequality becomes an equality. Therefore, $\mathrm{Slo}_\rho$ is of real dimension $2(g-1)\dim \mathfrak{g}^\textbf{R}$.
\end{proof}

Combining the results from above, we prove Theorem \ref{thm:characterization}:
\begin{proof}[Proof of Theorem \ref{thm:characterization}]
Given the argument in the paragraph directly before Theorem \ref{thm:characterization}, it suffices to prove the converse direction. Suppose $\mathrm{Slo}_\rho$ is a union of connected components, we show that $\hat{\rho}$ must be extended magical. If all of the connected components in $\mathrm{Slo}_\rho$ achieves expected dimension, Proposition \ref{prop:even_magical_triples} shows that $\hat{\rho}$ must be an even magical triple. Otherwise, there is a connected component in $\mathrm{Slo}_\rho$ that has strictly smaller dimension than the expected dimension. 

Under the nonablian Hodge correspondence, this implies that the surface group representations in this component are rigid in the sense of \cite{Kim_Pansu_2014}. From the result of Kim--Pansu in Remark \ref{rmk:Kim-Pansu} we deduce that, given our assumptions on the genus of the underlying curve and on the target Lie group $G^\textbf{R}$ to be the real form of a simple complex Lie group, the surface group representations in this component must be maximal representations into a Hermitian Lie group of nontube type. Applying the nonabelian Hodge correspondence again, this
implies that this connected component of $\mathrm{Slo}_\rho$ must be the maximal components of $\mathcal{M}(G^\textbf{R})$ for a nontube type Hermitian Lie group $G^\textbf{R}$.

We show that in the latter case, $\hat{\rho}$ must be an odd magical triple. Suppose not, then either $f+V_{even}\not\subset \mathfrak{m}$ or $\mathfrak{c}\not\subset \mathfrak{h}$, and therefore $\hat{\rho}$ \textit{cannot} be an even magical triple of $G_T^\textbf{R}$. By Proposition \ref{prop:part2}, we know that a maximal $G^\textbf{R}$-Higgs bundle can be written as $\mathcal{E}_{G^\textbf{R}}=\mathcal{E}_{G_t^\textbf{R}}\times_Z \mathcal{E}_{\hat{C}^\textbf{R}}$, where $\mathcal{E}_{G_t^\textbf{R}}$ has maximal Toledo number. Now $G_t^\textbf{R}$ is the real form of a simple complex Lie group, so we may apply our first argument and reason that since $\hat{\rho}$ does \textit{not} restrict to an even magical triple of $G_t^\textbf{R}$, there exists a maximal $\mathcal{E}_{G_t^\textbf{R}}\times_Z \mathcal{E}_{\hat{C}^\textbf{R}}$ that is not in the image of the restricted Slodowy map of $\rho$. Hence the Slodowy slice $\mathrm{Slo}_\rho \subset \mathcal{M}(G^\textbf{R})$ cannot be the union of maximal components, a contradiction.
\end{proof}

\section{A Cayley correspondence for nontube type $G^\textbf{R}$-Higgs bundles}
We conclude this paper by proving a Cayley correspondence for nontube type $G^\textbf{R}$-Higgs bundles. To the author's best knowledge, such a characterization was not described in any of the previous works (e.g. \cite{Bradlow_García-Prada_Gothen_2007} or \cite{Biquard_García-Prada_Rubio_2017}) concerning the Cayley correspondence for Hermitian Lie groups. Note that proving the geometric properties of injectivity, openness and closedness is a stronger statement than showing that the image of the restricted Slodowy map is a union of connected components. 

\begin{theorem}
\label{thm:cayley_correspondence_odd}
    Let $G$ a complex simple Lie group and let $\hat{\rho}$ be an odd magical triple of the real form $G^\textbf{R}\subset G$. Let $\tilde{G}^\textbf{R}$ be the subgroup of $G^\textbf{R}$ with Lie algebra $\tilde{\mathfrak{g}}^\textbf{R} \subset \mathfrak{g}^\textbf{R}$ and maximal compact subgroup $\mathcal{C}:=C\cap G^\textbf{R}$. Then the restriction of the Slodowy map to the preimage of $\mathcal{H}(G^\textbf{R})$ descends to an injective, open and closed map on the level of moduli spaces mapping onto the maximal components
    \begin{align*}
\Psi_\rho:\,\,&\mathcal{M}_{K^2}(\tilde{G}^\textbf{R})\times H^0(K) \to \mathcal{M}_{\mathrm{max}}(G^\textbf{R})\subset\mathcal{M}(G^\textbf{R}), \\
&((E_\mathcal{C},\varphi_\mathcal{C}), \,\varphi) \mapsto (E_H:=(E_\mathcal{C}\times E_T)[H], \underline{f}+\psi_\mathcal{C}+\varphi),
    \end{align*}
where $\psi_\mathcal{C}=\mathrm{ad}_f^{-2}(\phi_\mathcal{C})\in H^0(X,E_H[\mathfrak{m}]\otimes K)$. 
\end{theorem}

\begin{proof}
    The proof largely follows from the arguments in \cite{Bradlow_Collier_García-Prada_Gothen_Oliveira_2024} with minimal modifications. We explain here how their arguments apply to our situation. By Theorem \ref{thm:maximal_components}, we know that the image of $\Psi_\rho$ are the maximal components. We use the fact that $\mathcal{M}_{\mathrm{max}}(G^\textbf{R})\cong \mathcal{M}_{\mathrm{max}}(G_T^\textbf{R})$ and that any odd magical triple of $G^\textbf{R}$ restricts to an even magical triple of the tube type Hermitian Lie group $G_T^\textbf{R}$ (Proposition \ref{prop:restriction_of_odd_magical}). 
    
    Then injectivity of the Slodowy map on the level of Higgs pairs follows from Case (2) of Proposition 5.11 of \cite{Bradlow_Collier_García-Prada_Gothen_Oliveira_2024}, where we have $K$-twisting $m_c=2$ and the weights $\{l_j\}=\{0\}$. That the Slodowy map descends to an injective map on the level of moduli spaces follows from Theorem 7.7 of \cite{Bradlow_Collier_García-Prada_Gothen_Oliveira_2024}, which in turn follows from Lemma 7.9 of \cite{Bradlow_Collier_García-Prada_Gothen_Oliveira_2024}, again citing Case (2) for Hermitian Lie groups of tube type.
    
    To show that $\Psi_\rho$ is open, we apply the arguments in Section 7.3 of \cite{Bradlow_Collier_García-Prada_Gothen_Oliveira_2024}, which only uses the fact that an even magical triple of $\mathfrak{g}^\textbf{R}$ satisfies the conditions $\ker(\mathrm{ad}_f:\mathfrak{h}\to \mathfrak{m})=\mathfrak{c}$ and that $\mathrm{ad}_f:\mathrm{ad}_f(\mathfrak{m})\to \mathrm{ad}^2_f(\mathfrak{m})$ is an isomorphism. Finally, that $\Psi_\rho$ is closed follows from the argument in Section 7.4 of \cite{Bradlow_Collier_García-Prada_Gothen_Oliveira_2024} for the even magical triple of $G_T^\textbf{R}$.
\end{proof}

\printbibliography[heading=bibintoc, title=References]

\end{document}